\documentclass[11pt]{amsart}
\usepackage[utf8]{inputenc}
\usepackage[T1]{fontenc}
\usepackage{mathtools,amssymb}
\usepackage{tikz-cd}
\usepackage[capitalize,noabbrev]{cleveref}

\newtheorem{theorem}{Theorem}
\newtheorem{corollary}[theorem]{Corollary}
\newtheorem{lemma}[theorem]{Lemma}
\newtheorem{proposition}[theorem]{Proposition}
\newtheorem{question}[theorem]{Question}

\newtheorem{mainthm}{Main Theorem}
\newtheorem{fact}[theorem]{Fact}

\theoremstyle{definition}
\newtheorem{definition}[theorem]{Definition}

\newenvironment{subproof}[1][Proof]{\begin{proof}[#1]}{\end{proof}}

\newcommand{\st}{\, ;\,}
\newcommand{\set}[2]{\{#1 \st #2\}}
\newcommand{\seq}[2]{\langle #1 \st #2\rangle}
\newcommand{\leftexp}[2]{{\vphantom{#2}}^{#1}{#2}}
\newcommand{\funcs}[2]{\leftexp{#1}{#2}}
\newcommand{\forces}{\Vdash}
\newcommand{\rest}{\mathbin{\upharpoonright}}
\DeclareMathOperator{\CP}{\mathrm{CP}}
\DeclareMathOperator{\LCP}{\mathrm{LCP}}
\renewcommand{\P}{\mathbb{P}}
\newcommand{\Q}{\mathbb{Q}}
\newcommand{\R}{\mathbb{R}}
\renewcommand{\S}{\mathbb{S}}
\newcommand{\U}{\mathcal{U}}
\newcommand{\A}{\mathbb{A}}
\newcommand{\C}{\mathbb{C}}
\newcommand{\lexl}{<_{\mathrm{lex}}}
\DeclareMathOperator{\Term}{Term}

\DeclareMathOperator{\cf}{cf}

\DeclareMathOperator{\Add}{Add}

\DeclareMathOperator{\dom}{dom}

\DeclareMathOperator{\power}{\mathcal{P}}
\DeclareMathOperator{\Ult}{Ult}
\DeclareMathOperator{\supp}{supp}
\newcommand{\Ptail}{\P_{\textup{tail}}}
\newcommand{\Gtail}{G_{\textup{tail}}}

\begin{document}
	\title{Capturing sets of ordinals by normal ultrapowers}
	\author{Miha E.\ Habi\v{c}}
	\address[M.~E.~Habi\v{c}]{
		Bard College at Simon's Rock\\
		84 Alford Road\\
		Great Barrington, MA 01230\\
		USA
	}
	\email[Corresponding author]{mhabic@simons-rock.edu}
	\urladdr{https://mhabic.github.io}
	
	\author{Radek Honz{\'i}k}
        \address[R. Honz{\'i}k]{
		Department of Logic\\
		Faculty of Arts\\
		Charles University\\
		n\'am.\ Jana Palacha 2\\
		116 38 Praha 1\\
		Czech Republic
	}
	\email{radek.honzik@ff.cuni.cz}
	\urladdr{logika.ff.cuni.cz/radek}

	\thanks{The authors are grateful to Arthur Apter, Omer Ben-Neria, James Cummings, and Kaethe Minden for helpful discussions regarding the topics in this paper. We also thank the referee for their careful reading of the manuscript and their suggestions, which significantly improved the exposition in this paper.\\
		The first author was supported in part by the ESIF, EU Operational Programme Research, Development and Education, the International Mobility of Researchers in CTU project no.~(CZ.02.2.69/0.0/0.0/16\_027/0008465) at the Czech Technical University in Prague, and the joint FWF--GA\v{C}R grant no.~17-33849L: Filters, Ultrafilters and Connections with Forcing.\\ The second author was supported by GA{\v C}R-FWF grant \emph{Compactness principles and combinatorics}.\\
	}
	\begin{abstract}
		We investigate the extent to which ultrapowers by normal measures on \(\kappa\)
		can be correct about powersets \(\power(\lambda)\) for \(\lambda>\kappa\).
		We consider two versions of this question, the \emph{capturing property}
		\(\CP(\kappa,\lambda)\) and the \emph{local capturing property} \(\LCP(\kappa,\lambda)\). 
		Both of these describe the extent to which subsets of \(\lambda\) appear
		in ultrapowers by normal measures on \(\kappa\).
		After examining the basic properties of these two notions,
		we identify the exact consistency strength of \(\LCP(\kappa,\kappa^+)\).
		Building on results of Cummings, who determined the exact consistency strength
		of \(\CP(\kappa,\kappa^+)\), and using a variant of a forcing due to Apter and Shelah, we show
		that \(\CP(\kappa,\lambda)\) can hold at the least measurable cardinal.
	\end{abstract}
	\keywords{Capturing property, local capturing property, measurable cardinal, normal measure}
	\subjclass[2020]{03E55, 03E35}
	\maketitle
	
\section{Introduction}

It is well known that the ultrapower of the universe by a normal measure on some cardinal \(\kappa\) cannot
be very close to \(V\); for example, the measure itself never appears in the ultrapower. It follows that
these ultrapowers cannot compute \(V_{\kappa+2}\) correctly. In the presence of GCH, this is equivalent to
saying that the ultrapower is incorrect about \(\power(\kappa^+)\). But if GCH fails, it becomes conceivable
that a normal ultrapower could compute additional powersets correctly. This conjecture turns out to be
correct: Cummings~\cite{Cummings1993:StrongUltrapowersLongCoreModels}, answering a question of
Steel, showed that it is relatively consistent that there is a measurable cardinal
\(\kappa\) with a normal measure whose ultrapower computes \(\power(\kappa^+)\) correctly; in fact he showed
that this situation is equiconsistent with a \((\kappa+2)\)-strong cardinal \(\kappa\). 
In this paper we will study this \emph{capturing property} and its local variant further.

\begin{definition}
	Let \(\kappa\) and \(\lambda\) be cardinals. We say that the \emph{local capturing property}
	\(\LCP(\kappa,\lambda)\) holds if, for any \(x\subseteq\lambda\), there is a normal measure
	\(U_x\) on \(\kappa\) such that \(x\in \Ult(V,U_x)\). We shall say that
	\(U_x\) (or \(\Ult(V,U_x)\)) \emph{captures} \(x\).
\end{definition}

The full capturing property will amount to having a uniform witness for the local version.

\begin{definition}
	Let \(\kappa\) and \(\lambda\) be cardinals. We say that the \emph{capturing property}
	\(\CP(\kappa,\lambda)\) holds if there is a normal measure on \(\kappa\) that
	captures all subsets of \(\lambda\); in other words, a normal measure \(U\) such that
	\(\power(\lambda)\in \Ult(V,U)\).
\end{definition}

Some quick and easy observations: increasing \(\lambda\) clearly gives us stronger properties,
\(\CP(\kappa,\lambda)\) implies \(\LCP(\kappa,\lambda)\), and \(\CP(\kappa,\kappa)\)
holds for any measurable cardinal \(\kappa\).

Using this language, we can summarize Cummings' result as follows:

\begin{theorem}[Cummings]
	\label{thm:CummingsCP}
	If \(\kappa\) is \((\kappa+2)\)-strong, then there is a forcing extension in which
	\(\CP(\kappa,\kappa^+)\) holds. Conversely, if \(\CP(\kappa,\kappa^+)\) holds,
	then \(\kappa\) is \((\kappa+2)\)-strong in an inner model.
\end{theorem}

We should mention that \(\CP(\kappa,2^\kappa)\) is provably false: if it held, then some normal
ultrapower would contain all families of subsets of \(\kappa\), in particular the measure from which
it arose, which is impossible. Therefore a failure of GCH is necessary for \(\CP(\kappa,\kappa^+)\)
to hold. By work of Gitik~\cite{Gitik1989:NegationOfSCH}, this means that \(\CP(\kappa,\kappa^+)\)
has consistency strength at least that of a measurable cardinal \(\kappa\) with Mitchell rank \(o(\kappa)=\kappa^{++}\),
and the actual consistency strength of a \((\kappa+2)\)-strong cardinal \(\kappa\) is only slightly
beyond that.

The following are the main results of this paper. In \cref{sec:LCP} we analyse the consistency strength
of \(\LCP(\kappa,\kappa^+)\) and show that it is only a small step below the strength
of the full capturing property.

\begin{mainthm}
	\label{mainthm:mainThmLCP}
	Assuming GCH, if \(\LCP(\kappa,\kappa^+)\) holds, then \(o(\kappa)=\kappa^{++}\).
	Conversely, if \(o(\kappa)\geq\kappa^{++}\), then \(\LCP(\kappa,\kappa^+)\)
	holds in an inner model.
\end{mainthm}

In \cref{sec:CP} we continue the analysis in the case that GCH fails at \(\kappa\) and
show that the first part of the previous theorem, namely that \(\kappa\) has high Mitchell rank, fails dramatically if \(2^\kappa>\kappa^+\).

\begin{mainthm}
	\label{mainthm:mainThmCP}
	If \(\kappa\) is \((\kappa+2)\)-strong, then there is a forcing extension in which \(\CP(\kappa,\kappa^+)\)
	holds and \(\kappa\) is the least measurable cardinal.
\end{mainthm}

This last theorem is a nontrivial improvement of Cummings' result. Since the forcing he used to
achieve \(\CP(\kappa,\kappa^+)\) was relatively mild, \(\kappa\) remained quite large in the
resulting model; for example, it was still a measurable limit of measurable cardinals. Our theorem
shows that, while \(\CP(\kappa,\kappa^+)\) has nontrivial consistency strength, it does not directly
imply anything about the size of \(\kappa\) in \(V\) (beyond \(\kappa\) being measurable).

We will list questions that we have left open wherever appropriate throughout the paper.

\section{The local capturing property}
\label{sec:LCP}

Let us begin our analysis of the local capturing property with some simple observations.

\begin{lemma}
	\label{lemma:LCPCardinalAgreement}
	If \(\LCP(\kappa,\lambda)\) holds, then it can be witnessed by measures \(U\) for which
	\(\Ult(V,U)\) and \(V\) agree on cardinals up to and including \(\lambda\).
\end{lemma}

\begin{proof}
	Using a pairing function we can code a family of bijections \(f_\alpha\colon \alpha\to |\alpha|\)
	for \(\alpha\leq\lambda\) as a single subset \(y\subseteq \lambda\). 
	If we want to capture \(x\subseteq\lambda\) in an ultrapower as in the lemma, we simply
	capture (a disjoint union of) \(x\) and \(y\) using \(\LCP(\kappa,\lambda)\).
\end{proof}

\begin{proposition}
	\(\LCP(\kappa,(2^\kappa)^+)\) fails for any measurable \(\kappa\).
\end{proposition}

\begin{proof}
	If \(\LCP(\kappa,(2^\kappa)^+)\) held, there would have to be a normal measure ultrapower
	\(j\colon V\to M\) with critical point \(\kappa\) such that \(M\) was correct about cardinals
	up to and including \((2^\kappa)^+\), by \cref{lemma:LCPCardinalAgreement}. But no such
	ultrapower can exist, since the ordinals \(j(\kappa)\) and \(j(\kappa^+)\) are cardinals in \(M\)
	and both have size \(2^\kappa\) in \(V\).
\end{proof}

The following lemma is quite well known, but it will be key in many of our observations.

\begin{lemma}
	\label{lemma:factorMapCriticalPoint}
	Suppose that \(j\colon V\to M\) is an elementary embedding with critical point \(\kappa\)
	and consider the diagram
	\[\begin{tikzcd}
	V \arrow[r, "j"] \arrow[rd, "i", swap] & M\\ & N \arrow[u, "k", swap]
	\end{tikzcd}\]
	where \(i\) is the ultrapower by the normal measure on \(\kappa\)
	derived from \(j\) and \(k\) is the factor map. Then the critical point
	of \(k\) is strictly above \((2^\kappa)^N\).
\end{lemma}

\begin{proof}
	It is clear that the critical point of \(k\) is above \(\kappa\). Consider some ordinal \(\alpha\leq (2^\kappa)^N\).
	Fix a surjective map \(f\colon\power(\kappa)\to\alpha\) in \(N\) (and note that both \(N\) and \(M\) compute
	\(\power(\kappa)\) correctly). Since every ordinal up to and including \(\kappa\) is fixed by \(k\), it follows
	that \(k(f)=k\circ f\) is a surjection from \(\power(\kappa)\) to \(k(\alpha)\) and so \(k\rest\alpha\) is
	a surjection onto \(k(\alpha)\). It follows that we must have \(k(\alpha)=\alpha\).
\end{proof}

Using an old argument of Solovay, we can see that the optimal local capturing property automatically
holds at sufficiently large cardinals.

\begin{proposition}
	\label{prop:scHasMaxLCP}
	If a cardinal \(\kappa\) is \(2^\kappa\)-supercompact, witnessed by an embedding \(j\colon V\to M\),
	then \(\LCP(\kappa,2^\kappa)\) holds in both \(V\) and \(M\).
\end{proposition}

\begin{proof}
	We first show that \(\LCP(\kappa,2^\kappa)\) holds in \(V\). Suppose it fails. Then there is some
	\(x\subseteq 2^\kappa\) which is not captured by any normal measure on \(\kappa\). The model \(M\)
	agrees that this is the case, since it has all the normal measures on \(\kappa\) and all the functions
	\(f\colon\kappa\to\power(\kappa)\) that could represent \(x\). 
	Let \(i\) and \(k\) be as in \cref{lemma:factorMapCriticalPoint}.
	By that same lemma, the model \(N\) computes \(2^\kappa\) correctly
	and it also believes that there is some \(y\subseteq 2^\kappa\) which is not captured by any normal measure
	on \(\kappa\). This \(y\) is fixed by \(k\), so \(M\) also believes that \(y\) is not captured by any normal
	measure on \(\kappa\), and \(V\) agrees. But this is a contradiction, since \(y\) is captured by the ultrapower \(N\).
	Therefore \(\LCP(\kappa,2^\kappa)\) holds in \(V\).
	
	Observe that \(\LCP(\kappa,2^\kappa)\) only depends on \(\power(2^\kappa)\), the normal measures on \(\kappa\),
	and the representing functions \(\kappa\to\power(\kappa)\). The ultrapower \(M\) has all of these objects, therefore
	\(M\) must agree that \(\LCP(\kappa,2^\kappa)\) holds.
\end{proof}

In particular, if \(\kappa\) is \(2^\kappa\)-supercompact, then there are many \(\lambda<\kappa\) for which
\(\LCP(\lambda,2^\lambda)\) holds. 

The above argument seems to break down if \(\kappa\) is only \(\theta\)-supercompact for some \(\theta<2^\kappa\), even
if we are only aiming to capture subsets of \(\theta\); one simply cannot
conclude that \(M\) has all the necessary measures to correctly judge whether a set is a counterexample
to \(\LCP(\kappa,\lambda)\) or not. Thus, the following question remains open.

\begin{question}
	Suppose that \(\kappa\) is \(\theta\)-supercompact for some \(\kappa<\theta<2^\kappa\).
	Does it follow that \(\LCP(\kappa,\theta)\) holds?
\end{question}

The same conclusion as in \cref{prop:scHasMaxLCP} follows even if
\(\kappa\) is merely \((\kappa+2)\)-strong.

\begin{proposition}
	If a cardinal \(\kappa\) is \((\kappa+2)\)-strong, witnessed by an embedding \(j\colon V\to M\), then \(\LCP(\kappa,2^\kappa)\)
	holds in both \(V\) and \(M\).
\end{proposition}

\begin{proof}
	The argument works just like in \cref{prop:scHasMaxLCP}. 
	Note that \(M\) has all the functions
	\(\kappa\to\power(\kappa)\) and all the normal measures on \(\kappa\). Furthermore,
	\(M\) has all the subsets of \(2^\kappa\) (use a wellorder of \(V_{\kappa+1}\) in \(V_{\kappa+2}\)
	of ordertype \(2^\kappa\)). It follows that \(V\) and \(M\) have all the same counterexamples
	to \(\LCP(\kappa,2^\kappa)\).
\end{proof}

Reflecting back from \(M\) to \(V\), this last proposition implies that below a \((\kappa+2)\)-strong cardinal \(\kappa\) there are many cardinals \(\lambda\) satisfying
\(\LCP(\lambda,2^\lambda)\). This observation, together with Cummings' 
\cref{thm:CummingsCP}, tells us that the consistency strength of 
\(\LCP(\kappa,\kappa^+)\) is strictly lower than that of \(\CP(\kappa,\kappa^+)\). 
Let us determine this consistency strength exactly.

Recall that the \emph{Mitchell order} \(\triangleleft\) on a measurable cardinal \(\kappa\) is a relation on the
normal measures on \(\kappa\), where \(U\triangleleft U'\) if \(U\) appears in the ultrapower by \(U'\). It is a standard
fact that \(\triangleleft\) is wellfounded, and the \emph{Mitchell rank} of \(\kappa\) is the height \(o(\kappa)\) of this order.

\begin{proposition}
	\label{prop:LCPLargeMitchellRank}
	If \(\LCP(\kappa,2^\kappa)\) holds, then \(o(\kappa)=(2^\kappa)^+\).
\end{proposition}

\begin{proof}
	This is essentially the proof that the large cardinals mentioned in the previous two
	propositions have maximal Mitchell rank. We shall recursively build a Mitchell-increasing
	sequence \(\seq{U_\alpha}{\alpha<(2^\kappa)^+}\) of normal measures on \(\kappa\).
	So suppose that \(\seq{U_\alpha}{\alpha<\delta}\) has been constructed for some
	\(\delta<(2^\kappa)^+\). Using a pairing function we can code each measure \(U_\alpha\)
	as a subset of \(2^\kappa\), and then code the entire sequence \(\seq{U_\alpha}{\alpha<\delta}\)
	as a subset of \(2^\kappa\) as well. By \(\LCP(\kappa,2^\kappa)\) there is a normal measure
	\(U\) on \(\kappa\) which captures this subset, and thus the whole sequence of measures.
	We can then simply let \(U_\delta=U\).
\end{proof}

To show that the lower bound from this proposition is sharp we will pass to a suitable inner model.
Recall that a \emph{coherent sequence of normal measures} \(\U\) of length \(\lambda\) (where \(\lambda\) is an ordinal or \(\mathrm{Ord}\)) is given by a function
\(o^\U\colon \lambda\to\mathrm{Ord}\) and a sequence
\[\U=\seq{U_\alpha^\beta}{\alpha<\lambda, \beta<o^\U(\alpha)}\, ,\]
where each \(U_\alpha^\beta\) is a normal measure on \(\alpha\) and for each \(\alpha,\beta\),
if \(j_\alpha^\beta\) is the corresponding ultrapower map, we have
\[j_\alpha^\beta(\U)\rest \alpha+1 = \U\rest (\alpha,\beta)\, .\]
Here \(\U\rest (\alpha,\beta)= \seq{U_\gamma^\delta}{(\gamma,\delta)\lexl (\alpha,\beta)}\) and \(\U\rest\alpha = \U\rest (\alpha,0)\).

\begin{theorem}
	\label{thm:LUhasLCP}
	Suppose that \(V=L[\U]\) where \(\U\) is a coherent sequence of normal measures of length \(\kappa+1\)
	with \(o^\U(\kappa)=\kappa^{++}\). Then \(\LCP(\kappa,\kappa^+)\) holds.
\end{theorem}

\begin{proof}
	We shall show that, given any \(x\subseteq\kappa^+\), there is some \(\beta<\kappa^{++}\) such that
	\(x\in L[\U\rest (\kappa,\beta)]\). The theorem then immediately follows since, given \(x\),
	we can find a \(\beta\) as described, and the ultrapower by \(U_\kappa^\beta\) of \(L[\U]\)
	contains \(L[\U\rest (\kappa,\beta)]\), and therefore \(x\).
	
	So fix some \(x\subseteq\kappa^+\) and let \(\rho\) be a large regular cardinal so that \(x\in L_\rho[\U]\).
	Since GCH holds, we can find an elementary submodel \(M\prec L_\rho[\U]\) of size \(\kappa^+\)
	such that \(x,\U\in M\) and \(\kappa^+,\power(\kappa)\subseteq M\). Let \(\pi\colon M\to \overline{M}\) be
	the Mostowski collapse map.

	Note that \(\delta=M\cap\kappa^{++}=\pi(\kappa^{++})\) is an ordinal below \(\kappa^{++}\) and that all ordinals below \(\delta\) are fixed by \(\pi\). Moreover,
	\(\pi\) will fix all subsets of \(\kappa^+\) in \(M\) (since these can be described by sequences of ordinals of length \(<\delta\)), and therefore also
	all the measures \(U_\alpha^\beta\) for \((\alpha,\beta)\lexl (\kappa,\delta)\) (since each of these can be coded by a subset of \(\kappa^+\)). It follows that
	\(\pi(\U)\) is (in \(\overline{M}\)) a coherent sequence of normal measures of length \(\kappa+1\) 
	with
	\(o^{\pi(\U)}(\kappa)=\delta\), and that \(\pi(\U)=\U\rest(\kappa,\delta)\).
	Therefore
	\(\overline{M}=L_{\bar{\rho}}[\U\rest (\kappa,\delta)]\) for some \(\bar{\rho}<\rho\). Since
	\(x\subseteq \kappa^+\) was fixed by \(\pi\) as well, we get \(x\in \overline{M}\subseteq L[\U\rest (\kappa,\delta)]\).
\end{proof}

Even if, starting from a measurable cardinal \(\kappa\) of Mitchell order \(\kappa^{++}\), one could construct
a coherent sequence \(\U\) of normal measures with \(o^\U(\kappa)=\kappa^{++}\),
it seems to be an open question (according to~\cite{Steel2016:IntroToIteratedUltrapowers})
whether it is necessarily the case that \(\U\) remains coherent in \(L[\U]\). We avoid this issue
by using a result of Mitchell~\cite{Mitchell1983:SetsConstructedFromSequencesOfMeasuresRevisited},
who showed in ZFC that there is a sequence of filters \(\mathcal{F}\) (possibly empty, possibly of
length \(\mathrm{Ord}\), or anything in between) such that \(L[\mathcal{F}]\) satisfies GCH, \(\mathcal{F}\) is a coherent sequence of normal
measures in \(L[\mathcal{F}]\) and \(o^{\mathcal{F}}(\alpha)=\min (o(\alpha)^V, (\alpha^{++})^{L[\mathcal{F}]})\). The model we need will be exactly this \(L[\mathcal{F}]\).

\begin{corollary}
	\label{cor:MitchellRankGivesLCP}
	Assume that \(o(\kappa)\geq\kappa^{++}\). Then \(\LCP(\kappa,2^\kappa)\) holds in a transitive model of GCH.
\end{corollary}

\begin{proof}
	We may assume that \(\kappa\) is the largest measurable cardinal; if not,
	we can cut off the universe at the next inaccessible in order to achieve this.
	Let \(\mathcal{F}\) be the sequence of filters described above.
	By Mitchell's results we know that the sequence \(\mathcal{F}\) is a coherent sequence of normal measures in \(L[\mathcal{F}]\) and
	\(o^{\mathcal{F}}(\kappa)=(\kappa^{++})^{L[\mathcal{F}]}\).
	Since \(\kappa\) is the largest measurable, the length of
	\(\mathcal{F}\) is \(\kappa+1\), and it follows from \cref{thm:LUhasLCP} that
	\(\LCP(\kappa,\kappa^+)\) holds in \(L[\mathcal{F}]\).
%
\end{proof}

In fact, these canonical inner models satisfy a strong form of \(\LCP(\kappa,\kappa^+)\),
where there is a single function which represents any desired subset of \(\kappa^+\) in an
appropriate normal ultrapower.

\begin{definition}
	Let \(\kappa\) be a measurable cardinal. An \emph{\(H_{\kappa^{++}}\)-guessing Laver function} for \(\kappa\)
	is a function \(\ell\colon \kappa\to V_\kappa\) with the property that for any \(x\in H_{\kappa^{++}}\)
	there is an ultrapower embedding \(j\colon V\to M\) by a normal measure on \(\kappa\) such that \(j(\ell)(\kappa)=x\).
\end{definition}

It is obvious that the existence of an \(H_{\kappa^{++}}\)-guessing Laver function for \(\kappa\)
implies \(\LCP(\kappa,\kappa^+)\). The first author~\cite[Theorem 28]{Habic2019:JointDiamondsAndLaverDiamonds} showed
that this stronger property holds in appropriate extender models, in particular the one from 
\cref{cor:MitchellRankGivesLCP}.

Starting with a cardinal \(\kappa\) of high Mitchell rank, we obtained a model of the
local capturing property by passing to an inner model. We are unsure whether one can
obtain the local capturing property from the optimal hypothesis via forcing.

\begin{question}
	Suppose that GCH holds and \(o(\kappa)=\kappa^{++}\). Is there a forcing extension in which
	\(\LCP(\kappa,\kappa^+)\) holds?
\end{question}

It is important to note that the hypothesis in \cref{prop:LCPLargeMitchellRank}
is quite strong: we need to be able to capture all subsets of \(2^\kappa\) in order to be able
to conclude that the Mitchell rank of \(\kappa\) is large. One might wonder whether some large cardinal strength beyond measurability
can be derived even from weaker local capturing properties, for example \(\LCP(\kappa,\kappa^+)\)
assuming \(\kappa^+<2^\kappa\). As we shall see in the following section, the answer is
an emphatic no.

\section{The capturing property at the least measurable cardinal}
\label{sec:CP}

In this section we will give a proof of our second main theorem. Our argument owes a lot to
Cummings' original proof of \cref{thm:CummingsCP} and to the forcing machinery introduced
by Apter and Shelah. Nevertheless, we shall strive to give a mostly self-contained account, especially
with regard to the forcing notions used.

Let us first explain why we cannot simply use the proof from \cref{thm:CummingsCP}
and afterwards make \(\kappa\) into the least measurable cardinal
just by applying the standard methods of destroying measurable cardinals, such as iterated
Prikry forcing or adding nonreflecting stationary sets. In his argument, Cummings starts with a
\((\kappa,\kappa^{++})\)-extender embedding, lifts it through a certain iteration of Cohen forcings
(which will, among other things, ensure that \(2^\kappa>\kappa^+\), a necessary condition as we explained),
and concludes that the lifted embedding \(j\colon V[G]\to M[j(G)]\) is in fact equal to the ultrapower by some normal measure
on \(\kappa\) and \(M[j(G)]\) captures all the subsets of \(\kappa^+\) in the extension. One would now hope to be
able to lift this new embedding further, through any of the usual forcings which would make \(\kappa\)
into the least measurable cardinal. 
However, for this strategy to work, we should somehow ensure that \(\kappa\) is
not measurable in \(M[j(G)]\). Otherwise lifting the embedding through any of the usual forcing iterations to destroy all the measurables below \(\kappa\) over \(V[G]\) would require us to also destroy the measurability of \(\kappa\) over \(M[j(G)]\). But if we did that and maintained the capturing property at the same time, there would be enough agreement between the extensions of \(V[G]\) and \(M[j(G)]\) that \(\kappa\) would necessarily be nonmeasurable in the extension of \(V[G]\) as well. All this is to say that, since
\(\kappa\) is very much measurable in \(M[j(G)]\) after the forcing done by Cummings, a different 
approach is necessary.

Instead of first forcing the capturing property and then making \(\kappa\) into the least measurable,
the solution is to destroy all the measurable cardinals below \(\kappa\) and blow up \(2^\kappa\)
at the same time. The tools to make this approach work are due to Apter and 
Shelah~\cite{ApterShelah1997:MenasResultBestPossible, ApterShelah1997:OnStrongEqualityBetweenSupercompactnessStrongCompactness}.

\subsection{The forcing notions}

Let us review the particular forcing notions that will go into building our final forcing iteration.
The material in this subsection is contained, in some form or another, 
in Sections 1 of~\cite{ApterShelah1997:MenasResultBestPossible, ApterShelah1997:OnStrongEqualityBetweenSupercompactnessStrongCompactness}.

Since we will be discussing the strategic closure of some of these posets, let us fix some terminology.
If \(\P\) is a poset and \(\alpha\) is an ordinal, the \emph{closure game} for \(\P\) of length \(\alpha\)
consists of two players alternately playing conditions \(p\in\P\) in a descending sequence of length \(\alpha\),
with player II playing at limit steps. Player II loses the game if at any stage she is unable to make a move; otherwise
she wins. If \(\P\) is a poset and \(\kappa\) is a cardinal, we shall say that:

\begin{itemize}
	\item \(\P\) is \(\leq\kappa\)-strategically closed if player II has a winning strategy in the closure game
	for \(\P\) of length \(\kappa+1\).
	\item \(\P\) is \(\prec\kappa\)-strategically closed if player II has a winning strategy in the closure game
	for \(\P\) of length \(\kappa\).
	\item \(\P\) is \(<\kappa\)-strategically closed if it is \(\leq\lambda\)-strategically closed for all
	\(\lambda<\kappa\).
\end{itemize}

Recall that, if \(\kappa\) is a cardinal, a poset \(\P\) is called \emph{\(\leq\kappa\)-distributive} (or \emph{\(<\kappa\)-distributive}) if forcing with \(\P\) does not add any new sequences
of ordinals of length \(\leq\kappa\) (or \(<\kappa\)). This is equivalent to saying that
the intersection of \(\leq\kappa\) (or \(<\kappa\)) many open dense subsets of \(\P\)
is dense open.

If \(\kappa\) is a cardinal and \(\alpha\) is an ordinal, we let \(\Add(\kappa,\alpha)\)
be the usual forcing notion to add \(\alpha\) many Cohen subsets of \(\kappa\). We think
of conditions in \(\Add(\kappa,\alpha)\) as filling in a grid with \(\alpha\) many columns of height \(\kappa\) with 0s and 1s. Each condition is only allowed to fill in fewer than \(\kappa\) many
cells in the grid. Eventually, the generic will fill in the entire grid, and each column of
the grid will be a Cohen subset of \(\kappa\).

If \(\delta\geq\omega_2\) is a regular cardinal, we let \(\S_\delta\) be the forcing to add a nonreflecting
stationary subset of \(\delta\), consisting of points of countable cofinality.\footnote{In our argument we could
	use any other fixed cofinality below the large cardinal in question. We sacrifice a bit of generality in order to
	avoid carrying an extra parameter with us throughout the proof. The specific choice of countable cofinality also simplifies
	some arguments.}
A condition in \(\S_\delta\) is simply a bounded subset of \(x\subseteq\delta\), consisting of points of countable
cofinality and satisfying the property that \(x\cap\alpha\) is nonstationary in \(\alpha\) for every limit
\(\alpha<\delta\) of uncountable cofinality. The conditions in \(\S_\delta\) are ordered by end-extension. 
It is a standard fact that \(\S_\delta\) is
\(\prec\delta\)-strategically closed and, if \(2^{<\delta}=\delta\), is \(\delta^+\)-cc (see~\cite[Section 6]{Cummings:Handbook}
for more details). Note that the generic stationary set added will also be costationary, since it avoids
all ordinals of uncountable cofinality.

If \(S\subset\delta\) is a costationary set, let \(\C(S)\) be the forcing to shoot a club through \(\delta\setminus S\);
conditions are closed bounded subsets of \(\delta\setminus S\). Again, if \(2^{<\delta}=\delta\), then \(\C(S)\) will
be \(\delta^+\)-cc (\cite[Section 6]{Cummings:Handbook} has more details).
In the cases we will be interested in, \(\C(S)\) will also be \(<\delta\)-distributive (see \cref{lemma:NonRefStatShootClubIsCohen}).

Before we continue with the exposition, let us fix some terminology.

\begin{definition}
	Let \(\P\) and \(\Q\) be posets. We say that \(\P\) and \(\Q\) are \emph{forcing equivalent}
	if they have isomorphic dense subsets.
\end{definition}

This is not the most general definition of forcing equivalence that has appeared in the literature,
but it has the advantage of being obviously upward absolute between transitive models of set theory.

\begin{lemma}
	\label{lemma:NonRefStatShootClubIsCohen}
	If \(\delta\) is a cardinal satisfying \(\delta^{<\delta}=\delta\) then \(\S_\delta*\C(\dot{S})\), where
	\(\dot{S}\) is the name for the generic nonreflecting stationary set added by \(\S_\delta\), 
	is forcing equivalent to \(\Add(\delta,1)\).
\end{lemma}

\begin{proof}
	This is standard; the iteration has a dense \(<\delta\)-closed subset of size \(\delta\), which
	is equivalent to \(\Add(\delta,1)\) by~\cite[Theorem 14.1]{Cummings:Handbook}.
\end{proof}

Suppose that \(\gamma\) and \(\delta\) are regular cardinals, \(I\subseteq\delta\), and \(\vec{X}=\seq{x_\alpha}{\alpha\in I}\)
is a ladder system (meaning that each \(x_\alpha\subseteq\alpha\) is a \(\cf(\alpha)\)-sequence cofinal in \(\alpha\); the \(x_\alpha\) are called \emph{ladders}).
The forcing
\(\A(\gamma,\delta,\vec{X})\) consists of conditions \((p,Z)\) where

\begin{enumerate}
	\item \(p\) is a condition in the Cohen forcing \(\Add(\gamma,\delta)\), seen as filling in \(\delta\)
	many columns of height \(\gamma\) with 0s and 1s. We will denote by \(\supp(p)\subseteq\delta\) the set
	of indices of the nonempty columns of \(p\).
	\item \(p\) is a uniform condition, meaning that all of its nonempty columns have the same height.
	\item \(Z\) is a set of ladders from the ladder system \(\vec{X}\) and each ladder \(z\in Z\) is a subset of \(\supp(p)\).
\end{enumerate}

The conditions in \(\A(\gamma,\delta,\vec{X})\) are ordered by letting \((p',Z')\leq (p,Z)\)
if \(p'\leq p\) and \(Z'\supseteq Z\), and for any \(z\in Z\),
if \(\rho<\gamma\) is the index of a row that was empty in \(p\) but is nonempty in \(p'\), then both sets \(\{\iota\in z; p'(\iota,\rho)=0\}\) and \(\{\iota\in z; p'(\iota,\rho)=1\}\) are unbounded in \(\sup z\).
In other words, when strengthening the Cohen part of the condition, the \(z\in Z\) are promises that we will not add a row whose values stabilize when restricted to the columns indexed by \(z\).

The poset \(\A(\gamma,\delta,\vec{X})\) is similar to the poset \(P^1_{\delta,\lambda}[S]\) defined in~\cite[Section 4]{ApterShelah1997:MenasResultBestPossible},
with some differences which we believe will simplify the poset.
For example, our definition permits an arbitrary ladder system,
whereas Apter and Shelah work with a very specific one. For our applications, the specific
case studied by Apter and Shelah would have sufficed, but the poset can nevertheless be defined
more generally. We believe the additional generality will make the role of the side
conditions in the arguments more transparent and clarify where additional assumptions
on the parameters in the definition of \(\A(\gamma,\delta,\vec{X})\) are required.

Some comments are in order regarding the forcing \(\A(\gamma,\delta,\vec{X})\). It is similar
enough to the Cohen poset \(\Add(\gamma,\delta)\) that one would hope that it is just as simple to show
that this forcing also adds \(\delta\) new subsets of \(\gamma\) and so on. But with the addition
of the side conditions this is no longer clear. It is not even immediate that generically
we will fill out the entire binary matrix. On the other hand, if
we want to use this forcing as the main part of our construction to destroy many measurable cardinals,
then it cannot be too close to plain Cohen forcing after all.
This tension between the poset \(\A(\gamma,\delta,\vec{X})\) and the Cohen poset \(\Add(\gamma,\delta)\) is controlled
by the ladder system \(\vec{X}\), so we will have to choose these ladder systems carefully in our proof.

The following facts are parallel to the ones Apter and Shelah give in~\cite{ApterShelah1997:MenasResultBestPossible,ApterShelah1997:OnStrongEqualityBetweenSupercompactnessStrongCompactness}; we give proofs for the sake of completeness, but the reader familiar with their exposition should expect no surprises.

\begin{lemma}
	\label{lemma:ApterShelahCohenCompatibility}
	Two conditions \((p,Z)\) and \((q,W)\) in \(\A(\gamma,\delta,\vec{X})\) of equal height are compatible if and only if their Cohen parts are compatible.
\end{lemma}

\begin{proof}
	The forward implication is immediate. For the reverse, assume that \(p\) and \(q\) are compatible, so that \(r=p\cup q\) is a Cohen condition. Notice that
	\(r\) has the same height as \(p\) and \(q\), and that each nonempty column in \(r\) was either present already in both \(p\) and \(q\), or else it was present already in \(p\) and was empty in \(q\), or vice versa. If we let \(U=Z\cup W\), it then follows that \((r,U)\) is a common strengthening of both \((p,Z)\) and \((q,W)\). This is because, as far as ladders \(z\in Z\) are concerned, no new rows were added to the Cohen part when it was strengthened from \(p\) to \(r\), and similarly for \(W\).
\end{proof}

\begin{corollary}
	\label{cor:ApterShelahKnaster}
	For any regular \(\lambda>\gamma\), the poset \(\A(\gamma,\delta,\vec{X})\) is \(\lambda\)-Knaster (meaning that any collection of \(\lambda\) many conditions has a subset of size \(\lambda\) of pairwise compatible conditions)
	if and only if the poset \(\Add(\gamma,\delta)\) is \(\lambda\)-Knaster.
\end{corollary}

\begin{proof}
	For the forward direction, start with a family of Cohen conditions \(p_\alpha\) for \(\alpha<\lambda\), and associate to each the condition \((p_\alpha,\emptyset)\in \A(\gamma,\delta,\vec{X})\). Since \(\A(\gamma,\delta,\vec{X})\) is \(\lambda\)-Knaster, there is a subset \(J\subseteq \lambda\) of size \(\lambda\) such that the conditions \((p_\alpha,\emptyset)\) for \(\alpha\in I\) are pairwise compatible. This, of course, means that the conditions \(p_\alpha\) for \(\alpha\in I\)
	are pairwise compatible.
	
	Conversely, suppose that \(\Add(\gamma,\delta)\) is \(\lambda\)-Knaster and let \((p_\alpha,Z_\alpha)\) for \(\alpha<\lambda\) be conditions in \(\A(\gamma,\delta,\vec{X})\). If there are \(\lambda\) many conditions among these with the same Cohen part we are done, since all of those will be pairwise compatible. So let us assume that this doesn't happen. Since \(\lambda>\gamma\), there is a subset \(J\subseteq \lambda\) of size \(\lambda\) such that all conditions \((p_\alpha,Z_\alpha)\) have the same height.
	By our assumption and since \(\lambda\) is regular, we can thin out \(J\) further to assume that \(p_\alpha\neq p_\beta\) for distinct \(\alpha,\beta\in J\). Since \(\Add(\gamma,\delta)\) is \(\lambda\)-Knaster,
	we can thin out \(J\) even further until the conditions \(p_\alpha\) for
	\(\alpha\in J\) are all pairwise compatible. But since we already arranged them all to have the same height, \cref{lemma:ApterShelahCohenCompatibility} implies
	that the conditions \((p_\alpha,Z_\alpha)\) for \(\alpha\in J\) are also pairwise compatible in \(\A(\gamma,\delta,\vec{X})\).
\end{proof}

This corollary will be convenient when we need to gauge the chain condition or the Knasterness of the poset \(\A(\gamma,\delta,\vec{X})\). Typical applications will have \(\gamma\) inaccessible, \(\lambda\) a finite successor of \(\gamma\)
and GCH between \(\gamma\) and \(\lambda\). In those cases, a standard \(\Delta\)-system argument will guarantee that \(\Add(\gamma,\delta)\) (and therefore also \(\A(\gamma,\delta,\vec{X})\)) is \(\lambda\)-Knaster.

\begin{lemma}
	\label{lemma:ApterShelahClosed}
	Suppose that \(\gamma\) is regular, \(\delta\) is regular, and \(\vec{X}\)
	is a ladder system on some subset of \(\delta\). Then \(\A(\gamma,\delta,\vec{X})\)
	is \(<\gamma\)-closed.
\end{lemma}

The outright closure of the poset is a slight improvement over the presentation that
Apter and Shelah chose; they could only guarantee strategic closure, but the difference will not be significant.

\begin{proof}
	Start with a descending sequence of conditions \((p_\alpha,Z_\alpha)\) for \(\alpha<\lambda<\gamma\). We can get a candidate for a lower bound by simply taking
	unions in each coordinate, letting \(p=\bigcup_{\alpha<\lambda}p_\alpha\) and
	\(Z=\bigcup_{\alpha<\lambda}Z_\alpha\), but we need to verify that \((p,Z)\leq (p_\alpha,Z_\alpha)\). Consider any ladder \(z\in Z_\alpha\) and look at the restrictions \(p_\alpha\rest z\) and \(p\rest z\). For each new row in \(p\rest z\),
	we can find a \(\beta\) with \(\alpha<\beta<\lambda\) such that that row appears
	already in \(p_\beta\rest z\). But because we assumed that \((p_\beta,Z_\beta)\leq (p_\alpha,Z_\alpha)\), it must be the case that that row has unboundedly many 0s and 1s.
\end{proof}

Going forward, we will focus particularly on ladder systems supported on very sparse sets,
meaning those without any stationary initial segments. 
The following is essentially~\cite[Lemma 2]{ApterShelah1997:MenasResultBestPossible} and also~\cite[Lemma 2]{ApterShelah1997:OnStrongEqualityBetweenSupercompactnessStrongCompactness}: although
Apter and Shelah state the result for a very special ladder system, an inspection of their proof shows that the argument works in general. 

\begin{lemma}
	\label{lemma:ApterShelahDisjointify}
	Let \(\gamma\) be inaccessible and \(\delta\) regular. Suppose that \(I\subseteq\delta\) is nonstationary in
	its supremum and all of its initial segments are nonstationary in their suprema as well. Let \(\vec{X}\) be a ladder system on \(I\). Then there are (nonempty)
	final segments \(y_\alpha\) of each \(x_\alpha\in \vec{X}\) such that the \(y_\alpha\) are pairwise disjoint.
\end{lemma}

Let us briefly explain why this lemma will be useful. Suppose that we have a condition \((p,Z)\in\A(\gamma,\delta,\vec{X})\) and we would like to strengthen its Cohen part. We cannot just blindly extend \(p\), since the ladders in \(Z\) exert some control over what the rows of any extension of \(p\) might look like. If the ladders in \(Z\) are all pairwise disjoint, then this isn't a significant issue: we can fill in one cell of \(p\) and consider each ladder in \(Z\) separately,
filling in more of the row to ensure that unboundedly many 0s and 1s appear in the columns mentioned by \(z\). This naive strategy seems less solid when the ladders in \(Z\) overlap, since it might happen that, while satisfying the requirements given by one ladder \(z\in Z\), we inadvertently violate those given by a different \(z'\in Z\).

The point of \cref{lemma:ApterShelahDisjointify} is to allow us, under certain circumstances, to pretend that the ladders in \(Z\) really are pairwise disjoint. More precisely, instead of the ladders in \(z\in Z\), we are going to focus on their (pairwise disjoint) final segments \(y\) as provided by \cref{lemma:ApterShelahDisjointify}. The point is that, if we want to strengthen the Cohen part \(p\) of \((p,Z)\) to \(p'\) by filling in a cell in row \(\rho\), it suffices, for each \(z\in Z\), to have \(\{\iota\in y;p'(\iota,\rho)=0\}\) and \(\{\iota\in y; p'(\iota,\rho)=1\}\) be unbounded in \(\sup y=\sup z\), instead of the same requirement where \(y\) is replaced by \(z\), since \(y\) is a cofinal subset of \(z\).

\begin{lemma}
	\label{lemma:ApterShelahAddsReals}
	Suppose \(\gamma\) is inaccessible and \(\delta\geq \gamma\) is regular.
	Suppose that
	\(I\subseteq\delta\) and that all of its proper initial
	segments are nonstationary. Let \(\vec{X}\) be a ladder system on \(I\). Then a generic
	for \(\A(\gamma,\delta,\vec{X})\) is a total function on \(\delta\times\gamma\) and each of its
	columns is a new subset of \(\gamma\).
\end{lemma}

\begin{proof}
	We only need to show that, given
	a condition \((p,Z)\), we may extend that condition in order to fill any given empty 
	cell with an arbitrary value. 
	This is sometimes easy to do: if the height of \(p\) is \(\rho\) and we wish to fill a cell below height \(\rho\), we can simply fill that cell (and even its column up to height \(\rho\)) in whatever way we want. The reason is that the cell will only be empty if its whole column is empty (since \(p\) is uniform), but that means that no ladder in \(Z\) mentions that column. Consequently, the side condition plays no part when strengthening the Cohen part in that column.
	
	Let us now consider the case when we are attempting to fill in a cell in row \(\rho\), meaning the first new row above the height of \(p\).
	We start building
	the stronger Cohen condition by filling in that new cell in the desired way.
	We still need to make this new Cohen condition uniform (all previously nonempty columns now need to get an entry in row \(\rho\), but also the column we just added an entry to might have empty cells below row \(\rho\) if it was empty prior to this step), and pay attention to the
	promises we made regarding the ladders in \(Z\).
	
	If the column we just added to was empty, we can fill it up to height \(\rho\) in whatever way we want. The reason is the same as before: this column doesn't appear in \(\supp(p)\), so no ladder in \(Z\) mentions it, and therefore the side condition has nothing to say about how we extend this column. So let us focus on adding entries in row \(\rho\) to columns that were nonempty in \(p\).
	
	Since \(\gamma\) is inaccessible, \(Z\) has size less than \(\gamma\). Since
	\(\delta\geq\gamma\) is regular, the ladders in \(Z\) are bounded
	below \(\delta\) and we can pick some \(\delta'<\delta\) so that each ladder in \(Z\) is a subset of \(\delta'\). By assumption \(I\cap \delta'\) is nonstationary in \(\delta'\) and all of its initial segments are nonstationary in their suprema as well. It follows that we can apply \cref{lemma:ApterShelahDisjointify}
	to \(Z\) (seen as a ladder system on \(I\cap \delta'\)) in order to find pairwise disjoint final segments
	\(y\) of each \(z\in Z\).
	
	For each \(\alpha\in\supp(p)\), there is at most one such final segment \(y\) for which
	\(\alpha\in y\). If there is no such \(y\), we fill the cell in row \(\rho\) and column \(\alpha\) with a 0. On the other hand, if such a \(y\) exists,
	we fill the cells in row \(\rho\) and columns in \(y\) in an alternating pattern
	to make sure that there are unboundedly many 0s and 1s. The key fact is that these
	specifications do not contradict each other, since the sets \(y\) are pairwise disjoint. In this way, we extend \(p\) to a uniform condition \(p'\) of height \(\rho+1\).
	If we were filling in a cell in a previously nonempty column, then \(\supp(p')=\supp(p)\), and otherwise \(\supp(p')=\supp(p)\cup \{\beta\}\), where \(\beta\) is the index of the empty column in \(p\) that we filled cells in. Using the reasoning described after \cref{lemma:ApterShelahDisjointify}, it is also clear that \((p',Z)\leq (p,Z)\): given any \(z\in Z\), consider the values \(p'(\iota,\rho)\) for \(\iota\in z\). Let \(y\subseteq z\) be the associated final segment. Our construction made sure that \(p'(\iota,\rho)=0\) and \(p'(\iota',\rho)=1\)
	for unboundedly many \(\iota,\iota'\in y\), and therefore also for unboundedly many \(\iota,\iota'\in z\), so it follows that \((p',Z)\leq (p,Z)\).
	
	Having seen how to fill a cell in row \(\rho\) of a condition \((p,Z)\) of height \(\rho\), we can use the same process to fill a cell in any row above \(\rho\) as well. We simply use the same argument to increase the height of \(p\) one step at a time and the closure of \(\A(\gamma,\delta,\vec{X})\)
	from \cref{lemma:ApterShelahClosed} to pass through limit steps, until we reach the desired cell to be filled in.	
\end{proof}


If \(\delta\) is a regular cardinal and \(S\subseteq\delta\) is stationary, recall that a \(\clubsuit_\delta(S)\)-sequence
is a ladder system \(\seq{x_\alpha}{\alpha\in S}\) such that for any unbounded \(A\subseteq\delta\)
there is some \(\alpha\in S\) such that \(x_\alpha\subseteq A\).\footnote{The principle \(\clubsuit_\delta(S)\) is usually stated in the apparently stronger form where there are stationarily many \(\alpha\in S\) for which \(x_\alpha\subseteq A\). This formulation is equivalent to the one we use; see~\cite[Observation I.7.2]{Shelah1998:ProperImproperForcing}.}

\begin{lemma}
	\label{lemma:ApterShelahForcesNonmeasurable}
	Suppose that \(\gamma<\delta\) are regular cardinals, with \(\gamma\) inaccessible and 
	\(\theta^{<\gamma}<\delta\) for all \(\theta<\delta\).\footnote{Again, it is best to think of the case when \(\delta\) is a finite successor of \(\gamma\) and GCH holds for cardinals in the interval \([\gamma,\delta]\).}
	Let \(S\subseteq\delta\)
	be a nonreflecting stationary set consisting of points of countable cofinality, and let \(\vec{X}\)
	be a \(\clubsuit_\delta(S)\)-sequence. Then \(\A(\gamma,\delta,\vec{X})\) forces that \(\gamma\) is
	not measurable.
\end{lemma}

\begin{proof}
	The proof follows the strategy of~\cite[Lemma 3]{ApterShelah1997:OnStrongEqualityBetweenSupercompactnessStrongCompactness}. 
	We start with a condition \((p,Z)\)
	and a name \(\dot{U}\) for a countably complete ultrafilter on \(\gamma\). For each \(i<\delta\), let \(\dot{S}^1_i\) be the canonical name for the subset of \(\gamma\) whose characteristic function is given by the \(i\)th column of the generic, and let \(\dot{S}^0_i\) be the name for its complement. For each \(i<\delta\) we can find a stronger condition \((p_i,Z_i)\leq (p,Z)\) which decides whether \(\dot{S}^1_i\in \dot{U}\) or \(\dot{S}^0_i\in \dot{U}\). By thinning out if necessary, we may assume that all of the Cohen parts \(p_i\) have the same height, and that the value of \(k_i\in\{0,1\}\) determining which of \(\dot{S}^{k_i}_i\) is forced to be in \(\dot{U}\), is independent of \(i\) and equal to some \(k\in\{0,1\}\). Moreover, we may further strengthen these conditions to ensure that \(i\in\supp(p_i)\) for all \(i\).
	Our cardinal arithmetic assumption implies that \(\Add(\gamma,\delta)\) and, according to \cref{cor:ApterShelahKnaster}, \(\A(\gamma,\delta,\vec{X})\) are \(\delta\)-Knaster, so we can find
	an unbounded set \(I\subseteq\delta\) such that the conditions \((p_i,Z_i)\)
	for \(i\in I\) are pairwise compatible.
	
	We now use the \(\clubsuit_\delta(S)\)-sequence: there is an \(\alpha\in S\) for which
	\(x_\alpha\subseteq I\). We can now let \(p^*=\bigcup_{i\in x_\alpha} p_i\) and
	\(Z^*=\bigcup_{i\in x_\alpha} Z_i\).
	Since the conditions \(p_i\) were all pairwise compatible, \(p^*\) is also a Cohen condition, stronger than each \(p_i\). Moreover, since the conditions \(p_i\) and \(p^*\) have the same height, \((p^*,Z^*)\) is actually a strengthening of each \((p_i,Z_i)\) (the argument is the same as in the proof of \cref{lemma:ApterShelahCohenCompatibility}).
	
	Now consider the (even stronger) condition \((p^*,Z^*\cup\{x_\alpha\})\); it really is a condition since each \(i\in x_\alpha\) was in the support of \(p_i\), so \(x_\alpha\subseteq \supp(p^*)\). This
	condition forces that \(\dot{S}^k_i\in\dot{U}\) for \(i\in x_\alpha\). Since \(x_\alpha\) is countable, it follows that \((p^*,Z^*\cup\{x_\alpha\})\) also forces that the intersection \(\bigcap_{i\in x_\alpha}\dot{S}^k_i\) is in \(\dot{U}\).
	
	We now show that the condition \((p^*,Z^*\cup\{x_\alpha\})\) also forces that the above intersection is bounded in \(\gamma\). Suppose otherwise, that there is a stronger condition \((q,W)\) forcing that the intersection \(\bigcap_{i\in x_\alpha}\dot{S}^k_i\) has an element \(\xi\) above the height of \(p^*\). This can only happen if \(q\) has \(k\) as the entries in row \(\xi\) and columns \(i\in x_\alpha\). But this is impossible, since the definition of the ordering in \(\A(\gamma,\delta,\vec{X})\) and the fact that \((q,W)\leq (p^*,Z^*\cup\{x_\alpha\})\) require there to be unboundedly many 0 entries as well as unboundedly many 1 entries in row \(\xi\) and columns \(i\in x_\alpha\) in \(q\).
	
	This shows that the countably complete ultrafilter \(\dot{U}\) is forced to have a bounded element. Therefore it cannot be a normal ultrafilter on \(\gamma\), so \(\gamma\) is not measurable in the forcing	extension.
\end{proof}

%

The following lemma is~\cite[Lemma 1]{ApterShelah1997:MenasResultBestPossible} (and also~\cite[Lemma 1]{ApterShelah1997:OnStrongEqualityBetweenSupercompactnessStrongCompactness}); the reader may
find the proof there. The argument is much like the proof that \(\Add(\omega_1,1)\) forces
\(\diamondsuit\).

\begin{lemma}
	\label{lemma:NonRefForcesClub}
	Let \(\delta\) be a regular cardinal satisfying \(\delta^\omega=\delta\). Then \(\S_\delta\) forces
	that \(\clubsuit_\delta(S)\) holds, where \(S\) is the generic stationary set added.
\end{lemma}

Since we now know that \(\S_\delta\) adds a \(\clubsuit_\delta(S)\)-sequence, it makes sense to
consider the iteration \(\S_\delta*\A(\gamma,\delta,\vec{X})\), where \(\vec{X}\) is a
\(\clubsuit_\delta(S)\)-sequence added by the first stage of forcing. \cref{lemma:ApterShelahForcesNonmeasurable}
implies that this iteration will definitely make \(\gamma\) nonmeasurable (assuming we start from GCH or a similar hypothesis). The following lemma is a complement
to that result and can be used to show that the measurability of \(\gamma\) may be resurrected.
It corresponds to~\cite[Lemma 4]{ApterShelah1997:OnStrongEqualityBetweenSupercompactnessStrongCompactness}.

\begin{lemma}
	\label{lemma:iterationResolvesToCohen}
	Let \(\gamma<\delta\) be regular cardinals with \(\gamma\) inaccessible and \(\delta\)
	satisfying \(\delta^{<\delta}=\delta\). Then the iteration
	\(\S_\delta*(\A(\gamma,\delta,\vec{X})\times \C(\dot{S}))\), where \(\vec{X}\) is an arbitrary
	ladder system on \(S\), is equivalent to \(\Add(\delta,1)\times\Add(\gamma,\delta)\).
\end{lemma}

\begin{proof}
	We stick closely to the argument from~\cite{ApterShelah1997:OnStrongEqualityBetweenSupercompactnessStrongCompactness}.
	\cref{lemma:NonRefStatShootClubIsCohen} already told us that 
	\(\S_\delta *\C(\dot{S})\) is equivalent to \(\Add(\delta,1)\), so it only remains to
	show that, in the resulting extension \(V[S][C]\), \(\A(\gamma,\delta,\vec{X})^{V[S]}\) is equivalent to \(\Add(\gamma,\delta)^V=\Add(\gamma,\delta)^{V[S][C]}\).
	Since in \(V[S][C]\), the formerly stationary set \(S\) is no longer stationary,
	nor does it have any stationary initial segments, \cref{lemma:ApterShelahDisjointify} implies that we can disjointify the
	ladder system \(\vec{X}\) by picking final segments \(y_\alpha\subseteq x_\alpha\)
	for each \(\alpha\in S\).
	
	The set \(\delta\) can now be decomposed into the disjoint union of the \(y_\alpha\)
	plus the remainder \(R=\delta\setminus \bigcup_\alpha y_\alpha\). The key realization
	(as in the proof of \cref{lemma:ApterShelahAddsReals}) is that we can honour
	the promises given by a condition \((p,Z)\in \A(\gamma,\delta,\vec{X})^{V[S]}\)
	by strengthening \(p\) carefully on each \(y_\alpha\) (and these regions are pairwise disjoint and do not interfere with each other), and strengthening \(p\) quite freely
	on the remainder \(R\).
	
	To make this precise, let us write, given \(\alpha\in S\), \(\A_\alpha\) for the subposet of conditions \((p,\{y_\alpha\})\in \A(\gamma,\delta,\vec{X})^{V[S]}\) for which \(\supp(p)=y_\alpha\). Let us also write \(\A_R\) for the subposet of those conditions \((p,\emptyset)\) for which \(\supp(p)\subseteq R\).
	Each of the posets \(\A_\alpha\) is a \(<\gamma\)-closed poset with the induced ordering, and each has size \(\gamma\). This means that each \(\A_\alpha\) is equivalent to \(\Add(\gamma,1)\) by \cite[Theorem 14.1]{Cummings:Handbook}. On the other hand, since the conditions in \(\A_R\)
	have empty side conditions, the induced ordering there behaves exactly like
	\(\Add(\gamma,R)\). 
	
	Given a condition \((p,Z)\in \A(\gamma,\delta,\vec{X})^{V[S]}\), we can decompose it into
	the sequence of restrictions \((p\rest y_\alpha,\{y_\alpha\})\) and \((p\rest R,\emptyset)\). We would like to say that this decomposition gives rise to
	an isomorphism between \(\A(\gamma,\delta,\vec{X})^{V[S]}\) and \(\prod_{\alpha\in S}\A_\alpha\times \A_R\) (where the product is taken with \(<\gamma\)-support). Unfortunately, that is not quite the case: this map is not surjective, as its range consists exactly of those conditions in the product whose Cohen components in each factor have the same height.
	However, the range is still dense in the product, which shows that \(\A(\gamma,\delta,\vec{X})^{V[S]}\) and \(\prod_{\alpha\in S}\A_\alpha\times \A_R\)
	are equivalent.
	
	Putting together the equivalences from the last two paragraphs, we obtain an equivalence between
	\(\A(\gamma,\delta,\vec{X})^{V[S]}\) and \(\prod_{\alpha\in S}\Add(\gamma,1)\times \Add(\gamma,R)\), where the product \(\prod_{\alpha\in S}\Add(\gamma,1)\) is taken with \(<\gamma\)-support, and we can conclude that \(\A(\gamma,\delta,\vec{X})^{V[S]}\) is equivalent to \(\Add(\gamma,\delta)\).
\end{proof}

\subsection{Some additional facts about forcing and elementary embeddings}

In this subsection we collect some facts about forcing and ultrapowers, some more standard than others, that we will need throughout our paper. We indicate at each the 
parallel result from~\cite{Cummings1992:GCHAtSuccessorsFailsAtLimits}
or~\cite{Cummings:Handbook}, where proofs are also given.

\begin{fact}[{\cite[Proposition 9.1]{Cummings:Handbook}}]
	\label{fact:LiftingCriterion}
	Suppose that \(M\) and \(N\) are transitive models of ZFC and \(j\colon M\to N\) is an
	elementary embedding. Let \(\P\in M\) be a poset, let \(G\) be \(\P\)-generic over \(M\)
	and let \(H\) be \(j(\P)\)-generic over \(N\). If \(j[G]\subseteq H\) then \(j\) can be
	extended to an elementary embedding \(j\colon M[G]\to N[H]\) satisfying \(j(G)=H\).
\end{fact}

\begin{fact}[{\cite[Section 1.2.2, Fact 3]{Cummings1992:GCHAtSuccessorsFailsAtLimits}}]
	\label{fact:ExtenderLift}
	With the notation of the previous fact, if \(j\colon M\to N\) is a \((\kappa,\lambda)\)-extender embedding, then so is the lift \(j\colon M[G]\to N[H]\).
	In particular, if \(i\) is the ultrapower by a normal measure on \(\kappa\), then so
	is \(j\).
\end{fact}

\begin{fact}[{\cite[Section 1.2.2, Fact 2]{Cummings1992:GCHAtSuccessorsFailsAtLimits}}]
	\label{fact:Transfer}
	With the notation of the previous fact, suppose that \(j\) is a
	\((\kappa,\lambda)\)-extender embedding and that \(\P\) is \(\leq\kappa\)-distributive
	in \(M\). Then \(j[G]\) generates a \(j(\P)\)-generic filter over \(N\).
\end{fact}

\begin{fact}[{\cite[Proposition 8.1]{Cummings:Handbook}}]
	\label{fact:DiagCriterion}
	Let \(M\) be an inner model of ZFC, let \(\P\in M\) be a poset and let \(\kappa\) be a cardinal. 
	Suppose that \(\P\) is \(\prec\kappa\)-strategically closed (in \(V\)) and that the set of maximal antichains
	of \(\P\) in \(M\) has cardinality at most \(\kappa\) in \(V\). Then there is a \(\P\)-generic
	filter \(G\) over \(M\) in \(V\).
\end{fact}

\begin{fact}[{\cite[Section 1.2.3, Fact 3]{Cummings1992:GCHAtSuccessorsFailsAtLimits}}]
	\label{fact:ClosurePreservation}
	Let \(M\) be an inner model of ZFC, let \(\P\in M\) be a poset and let \(\kappa\) be a cardinal.
	Suppose that \(\funcs{\kappa}{M}\subseteq M\) and that \(\P\) is \(\kappa^+\)-cc in \(V\).
	Let \(G\) be \(\P\)-generic over \(V\). Then \(M[G]\) is an inner model of \(V[G]\)
	and \(\funcs{\kappa}{M[G]}\subseteq M[G]\) in \(V[G]\).
\end{fact}

Recall that if \(\P\) is a poset and \(\dot{\Q}\) is a \(\P\)-name for a poset, the
\emph{term forcing poset} \(\Term(\P,\dot{\Q})\) consists of \(\P\)-names for elements of \(\dot{\Q}\), ordered
by letting \(\sigma\leq\tau\) if \(\P\forces \sigma\leq\tau\). 

\begin{fact}[{\cite[Section 1.2.5, Fact 1]{Cummings1992:GCHAtSuccessorsFailsAtLimits}}]
	\label{fact:TermForcing}
	If \(G\subseteq\P\) and
	\(H\subseteq\Term(\P,\dot{\Q})\) are generic over \(V\), then
	\(\set{\sigma^G}{\sigma\in H}\subseteq \dot{\Q}^G\) is generic over \(V[G]\).
\end{fact}

\begin{lemma}[{\cite[Section 1.2.5, Fact 2]{Cummings1992:GCHAtSuccessorsFailsAtLimits}}]
	\label{lemma:TermForcingCohen}
	Suppose that \(\kappa\) is a cardinal satisfying \(\kappa^{<\kappa}=\kappa\) and let \(\P\) be a
	\(\kappa\)-cc forcing of size \(\kappa\). Let \(\dot{\Q}_\lambda\) be the \(\P\)-name for \(\Add(\kappa,\lambda)\)
	in the extension. Then \(\Term(\P,\dot{\Q}_\lambda)\) is forcing equivalent, in \(V\), to \(\Add(\kappa,\lambda)\).
\end{lemma}

\begin{lemma}
	\label{lemma:GitikMerimovich}
	Let \(\kappa\) be a measurable cardinal satisfying \(2^\kappa=\kappa^+\) and let \(j\colon V\to M\) be the ultrapower by a normal measure on \(\kappa\). Given any finite \(n\geq 1\), the forcings \(j(\Add(\kappa,\kappa^{+n}))\) and \(\Add(\kappa^+,\kappa^{+n})\)
	are equivalent in \(V\).
\end{lemma}

Cummings gave a proof of this lemma for \(n=2\) in~\cite[Section 1.2.6, Fact 2]{Cummings1992:GCHAtSuccessorsFailsAtLimits} (attributing the proof to Woodin),
and Gitik and Merimovich proved the generalization to all \(n\) in~\cite[Lemma 3.2]{GitikMerimovich1997:PossibleValuesContinuum}.

\begin{lemma}
	\label{lemma:ccDistributive}
	Let \(\kappa\) be a regular cardinal, let \(\P\) be some \(<\kappa\)-distributive forcing notion, and let \(\Q\) be a \(\kappa\)-cc forcing notion. 
	If \(\P\) forces that \(\Q\) is \(\kappa\)-cc, then
	\(\Q\) forces that \(\P\) is \(<\kappa\)-distributive.
\end{lemma}



\begin{proof}
	Let \(G\times H\) be \(\P\times\Q\)-generic over \(V\) and consider
	a sequence \(\vec{x}\) of ordinals in \(V[H][G]\) of some length less than \(\kappa\). We wish to see that \(\vec{x}\in V[H]\). Since \(\vec{x}\in V[H][G]=V[G][H]\) and \(\Q\) is \(\kappa\)-cc in \(V[G]\), we can find a nice \(\Q\)-name \(\sigma\) for \(\vec{x}\) in \(V[G]\) that can also be coded by a
	sequence of ordinals of length \(<\kappa\). Since \(\P\) is \(<\kappa\)-distributive, this name \(\sigma\) is already in \(V\), and so
	\(\vec{x}\) must appear in \(V[H]\), as desired.
\end{proof}

The following key observation was already implicit in Cummings' proof of \cref{thm:CummingsCP}.
It shows that, as long as one can arrange the value of \(2^\kappa\) appropriately, the apparently difficult part
of the capturing property tends to follow for free from the construction.

\begin{lemma}
	\label{lemma:extenderLiftsToUltrapower}
	Suppose that \(j\colon V\to M\) is a \((\kappa,\lambda)\)-extender embedding and \(2^\kappa\geq\lambda\).
	Then \(j\) is the ultrapower by a normal measure on \(\kappa\).
\end{lemma}

\begin{proof}
	Let \(i\colon V\to N\) be the ultrapower by the normal measure derived from \(j\) and let
	\(k\colon N\to M\) be the factor embedding. Consider some \(x\in M\). Since \(j\) is a
	\((\kappa,\lambda)\)-extender embedding, we can write \(x=j(f)(\alpha)\) for some \(\alpha<\lambda\)
	and some function with domain \(\kappa\). By \cref{lemma:factorMapCriticalPoint} the critical
	point of \(k\) is above \(\lambda\) and therefore
	\[x=j(f)(\alpha)=k(i(f))(\alpha)=k(i(f)(\alpha))\, ,\]
	which shows that \(k\) is surjective. On the other hand, \(k\) is an elementary embedding, so it is also injective. It follows that \(k\) is an isomorphism of transitive structures
	and thus trivial, so we can conclude that \(j=i\).
\end{proof}

\subsection{The proof}

We are now ready to prove the second main theorem. We restate it here for convenience.

\begin{theorem}
	\label{thm:mainThmCP}
	If \(\kappa\) is \((\kappa+2)\)-strong, then there is a forcing extension in which
	\(\CP(\kappa,\kappa^+)\) holds, \(2^\kappa=\kappa^{++}\), and \(\kappa\) is the least measurable.
\end{theorem}

This theorem shows that the hypothesis in \cref{prop:LCPLargeMitchellRank} is in some
sense optimal: if \(2^\kappa>\kappa^+\) then \(\LCP(\kappa,\kappa^+)\) 
is not enough to conclude that the Mitchell rank of \(\kappa\) is large. In fact, even 
\(\CP(\kappa,\kappa^+)\) can hold at the least measurable cardinal.

\begin{proof}
	We make some simplifying assumptions to start with. We may assume that GCH holds and that
	the \((\kappa+2)\)-strongness of \(\kappa\) is witnessed by a \((\kappa,\kappa^{++})\)-extender embedding
	\(j\colon V\to M\). We have the usual diagram
	\[\begin{tikzcd}
	V \arrow[r, "j"] \arrow[rd, "i", swap] & M\\ & N \arrow[u, "k", swap]
	\end{tikzcd}\]	
	where \(i\) is the induced normal ultrapower map. Using the GCH and \cref{lemma:factorMapCriticalPoint}, we can see that
	the critical point of \(k\) is \((\kappa^{++})^N\). Using the argument from~\cite{Cummings1993:StrongUltrapowersLongCoreModels},
	we may also assume that, in \(V\), there is an \(i(\Add(\kappa,\kappa^{++}))\)-generic filter over \(N\).
	
	The following observation will be important, and we include the straightforward proof.
	
	\begin{lemma}
		The map \(k\) is a \(((\kappa^{++})^N,\kappa^{++})\)-extender embedding. That is,
		\[M=\{k(g)(\alpha);\alpha<\kappa^{++}, \dom(g)=(\kappa^{++})^N\}\,.\]
	\end{lemma}
	
	\begin{subproof}
		We assumed that we could write \(M\) in the form \[M=\{j(f)(\alpha);\alpha<\kappa^{++}, \dom(f)=\kappa\}\,.\]
		Now take an arbitrary element \(j(f)(\alpha)\) of \(M\). We can rewrite it as \((k(i(f))(\alpha))\). If we
		now take \(g=i(f)\rest(\kappa^{++})^N\), it is not hard to see that \(j(f)(\alpha)=k(g)(\alpha)\), showing
		inclusion in one direction.
		
		For the other direction, take an element of the form \(k(g)(\alpha)\). The function \(g\) itself is of the
		form \(i(F)(\kappa)\) for some function \(F\) with domain \(\kappa\), since \(N\) is the ultrapower of \(V\) by
		a normal measure on \(\kappa\). This means we can write \(k(g)(\alpha)=k(i(F)(\kappa))(\alpha)= (j(F)(\kappa))(\alpha)\), since the critical point of \(k\) is above \(\kappa\). 
		Let \(\langle \cdot,\cdot\rangle\colon \kappa^{++}\times \kappa^{++}\to\kappa^{++}\) be a bijection, and define a function
		\(f\colon \kappa^{++} \to V\) by \(f(\xi)=F(\xi_1)(\xi_2)\), where \(\xi=\langle \xi_1,\xi_2\rangle\) and where the definition only makes sense if \(F(\xi_1)\) is a function with \(\xi_2\) in its domain (in other cases we can define \(f\) arbitrarily). It is now straightforward
		to see, using elementarity, that \(j(f)(\langle \kappa,\alpha\rangle)=(j(F)(\kappa))(\alpha)=k(g)(\alpha)\).
	\end{subproof}
	
	We now specify the forcing we will use. Let \(\P_\kappa\) be the Easton support iteration of length \(\kappa\)
	which forces at inaccessible \(\gamma<\kappa\) with \(\S_{\gamma^{++}}*\A(\gamma,\gamma^{++},\vec{X})\), where
	\(\vec{X}\) is some \(\clubsuit_{\gamma^{++}}(S)\)-sequence added by \(\S_{\gamma^{++}}\).\footnote{It does not matter much
	how we pick these \(\clubsuit\)-sequences. One possible way is to fix in advance a wellordering of some large \(H_\theta\)
	and always pick the least appropriate name.}
	Let \(G_\kappa\) be \(\P_\kappa\)-generic over \(V\). 
	We can factor \(j(\P_\kappa)\) as
	\[j(\P_\kappa)=\P_\kappa *\S_{\kappa^{++}}*\A(\kappa,\kappa^{++},\vec{Y})*\Ptail\,,\]
	where \(\vec{Y}\) is the \(\clubsuit_{\kappa^{++}}\)-sequence used by the forcing at stage
	\(\kappa\) in \(M[G_\kappa]\) and \(\Ptail\) is the remainder of the forcing between \(\kappa\) and \(j(\kappa)\).
	The full forcing that will give us our result is then
	\[\P=\P_\kappa*\S_{\kappa^{++}}*(\A(\kappa,\kappa^{++},\vec{Y})\times \C(\dot{S}))\,.\]
	Let us carefully try to lift the embedding \(j\) through this forcing.
	
	First, we can rewrite \(i(\P_\kappa)\) as
	\[i(\P_\kappa)=\P_\kappa*(\S_{\kappa^{++}}*\A(\kappa,\kappa^{++},\vec{Y}'))^{N^{\P_\kappa}}*\Ptail'\,,\]
	where \(\vec{Y}'\) and \(\Ptail'\) are defined similarly to \(\vec{Y}\) and \(\Ptail\) in the case of \(j(\P_\kappa)\) above.
	Since \(G_\kappa\) is generic over all of \(V\), it is definitely generic over \(N\) and \(M\). 
	The forcing \(\P_\kappa\) is below the critical point of the embedding \(k\), so we can easily lift 
	it to \(k\colon N[G_\kappa]\to M[G_\kappa]\). Moreover, since \(\P_\kappa\) is \(\kappa\)-cc, \(N[G_\kappa]\) will
	be closed under \(\kappa\)-sequences in \(V[G_\kappa]\).
	
	We now claim that, in \(V[G_\kappa]\), there is an \((\S_{\kappa^{++}})^{N[G_\kappa]}\)-generic over \(N[G_\kappa]\),
	and moreover that this generic amounts to a nonstationary subset of \((\kappa^{++})^N\) (which is an ordinal of cofinality \(\kappa^+\) in \(V\)) in \(V[G_\kappa]\). This follows from
	\cref{lemma:NonRefStatShootClubIsCohen}, which tells us that the iteration \(\S_{\kappa^{++}}*\C(\dot{S})\)
	is equivalent to \(\Add(\kappa^{++},1)\). Since \(V[G_\kappa]\) has an \(\Add(\kappa^{++},1)^{N[G_\kappa]}\)-generic
	over \(N[G_\kappa]\) (as this forcing is \(\leq\kappa\)-closed in \(V[G_\kappa]\) and only has \(\kappa^+\) many
	dense subsets from \(N[G_\kappa]\)), we can also extract the generic for \(\S_{\kappa^{++}}^{N[G_\kappa]}\).
	Furthermore, this generic stationary set will be nonstationary in \(V[G_\kappa]\), as witnessed by the generic club added by \(\C(\dot{S})\).
	
	So let \(S'\in V[G_\kappa]\) be \((\S_{\kappa^{++}})^{N[G_\kappa]}\)-generic over \(N[G_\kappa]\). 
	This means that \(S'\) is, in \(N[G_\kappa][S']\), a nonreflecting
	stationary subset of \((\kappa^{++})^{N[G_\kappa]}\). In particular, none of its
	proper initial segments are stationary in their supremum. This statement is upwards
	absolute, so \(V[G_\kappa]\supseteq N[G_\kappa][S']\) agrees about the nonstationarity
	of the initial segments of \(S'\). But more than this, \(S'\) itself is nonstationary
	in its supremum \((\kappa^{++})^N\), as we noted in the previous paragraph.
	Finally, observe that \((\kappa^{++})^N<\kappa^{++}\); this is because \(i(\kappa)>(\kappa^{++})^N\) and \(i(\kappa)\) has size \(2^\kappa=\kappa^+\) in \(V\),
	since \(i\) is the ultrapower by a normal measure on \(\kappa\).
	Together, these facts imply that \(S'\)
	is a condition in the real \(\S_{\kappa^{++}}\). Let \(S\) be some
	\(\S_{\kappa^{++}}\)-generic over \(V[G_\kappa]\) containing \(S'\). 
	The embedding \(k\) lifts again to
	\(k\colon N[G_\kappa][S']\to M[G_\kappa][S]\); this is because the critical point
	of \(k\) is \((\kappa^{++})^N\), which means that \(k[S']=S'\subseteq S\) by the
	choice of \(S\).
	
	Now consider the \(\clubsuit_{(\kappa^{++})^{N}}\)-sequence \(\vec{Y}'\) used by \(i(\P_\kappa)\) at
	stage \(\kappa\). Since the critical point of \(k\) is \((\kappa^{++})^N\), the sequence \(\vec{Y}'\) is
	simply an initial segment of the sequence \(\vec{Y}=k(\vec{Y}')\) used by \(j(\P_\kappa)\) at stage \(\kappa\).\footnote{We could have arranged matters so that \(\vec{Y}\) was also
	a \(\clubsuit_{\kappa^{++}}(S)\)-sequence in \(V[G_\kappa][S]\), but this will not be important
	for the argument.}
	It follows that, if we look at the forcing \(\A(\kappa,\kappa^{++},\vec{Y})\) in \(V[G_\kappa][S]\),
	we can write it as a product
	\begin{equation}
	\label{eq:factorization}
		\A(\kappa,\kappa^{++},\vec{Y})\cong \A(\kappa,(\kappa^{++})^N,\vec{Y}') \times \A(\kappa,\kappa^{++}\setminus (\kappa^{++})^N, \vec{Y}\setminus \vec{Y}')\,.
	\end{equation}
	There is a slight abuse of notation in the second factor, since the set
	 \(\kappa^{++}\setminus(\kappa^{++})^N\) is not an ordinal.
	Nevertheless, we trust that our meaning is clear.
	Observe also that, since \(\S_{\kappa^{++}}\) does not
	add bounded subsets to \(\kappa^{++}\), we know
	\[\A(\kappa,(\kappa^{++})^N,\vec{Y}')^{V[G_\kappa][S]} = \A(\kappa,(\kappa^{++})^N,\vec{Y}')^{V[G_\kappa]} = \A(\kappa,(\kappa^{++})^N,\vec{Y}')^{N[G_\kappa][S']}\, .\]

	Let \(g'\) be \(\A(\kappa,(\kappa^{++})^N,\vec{Y}')\)-generic over \(V[G_\kappa][S]\); in particular, it is also generic over
	\(N[G_\kappa][S']\). Since \(g'\) is generic for a forcing that is \(\kappa^{+}\)-cc in \(V[G_\kappa]\),
	it follows that \(N[G_\kappa][S'][g']\) is still closed under \(\kappa\)-sequences in \(V[G_\kappa][g']\).
	We can conclude from this that \(\Ptail'\) is \(\prec\kappa^+\)-strategically closed in \(V[G_\kappa][g']\). This is because this poset is such in \(N[G_\kappa][S'][g']\), being an Easton support iteration all of whose iterands are at least \(\prec \kappa^+\)-strategically closed according to \cref{lemma:ApterShelahClosed}. As we mentioned, the model \(N[G_\kappa][S'][g']\) is closed under \(\kappa\)-sequences in \(V[G_\kappa][g']\), and therefore the winning strategy in the closure game for \(\Ptail'\) of length \(\kappa^+\) in \(N[G_\kappa][S'][g']\) remains winning in the larger model \(V[G_\kappa][g']\) (since any losing play would be of length shorter than \(\kappa^+\) and available in the smaller model).
	This, together with the fact that \(\Ptail'\) has only \(\kappa^+\) many dense open subsets from \(V[G_\kappa][g']\)
	(and therefore only \(\kappa^+\) many maximal antichains), allows us to build, using \cref{fact:DiagCriterion} in \(V[G_\kappa][g']\), a \(\Ptail'\)-generic \(\Gtail'\) over
	\(N[G_\kappa][S'][g']\) and lift the embedding \(i\) to
	\[i\colon V[G_\kappa]\to N[G_\kappa][S'][g'][\Gtail']\, .\]

	We can now force over \(V[G_\kappa][S]\), using the factorization~\eqref{eq:factorization}, to complete \(g'\) to \(g\)
	which is fully \(\A(\kappa,\kappa^{++},\vec{Y})\)-generic over \(V[G_\kappa][S]\). In the extension \(V[G_\kappa][S][g]\)
	we can finally also lift the map \(k\) through the last two stages of forcing and obtain
	\[k\colon N[G_\kappa][S'][g'][\Gtail']\to M[G_\kappa][S][g][\Gtail]\, ,\]
	where \(\Gtail\) is the filter generated by the pointwise image of \(\Gtail'\).
	The lift through \(g'\) is straightforward: the critical point of \(k\) is
	\((\kappa^{++})^N\), so \(k[g']=g'\subseteq g\). On the other hand, the forcing \(\Ptail'\) is at least
	\(\leq(\kappa^{++})^N\)-strategically closed in \(N[G_\kappa][S'][g']\), so \cref{fact:Transfer} together
	with the knowledge that \(k\) is a \(((\kappa^{++})^N,\kappa^{++})\)-extender embedding show that the pointwise image
	of \(\Gtail'\) really does generate a generic filter.
	
	Composing the two lifts of \(i\) and \(k\)
	gives us a lift of \(j\). The situation is summarized in the following diagram; we should keep in mind that the pictured
	embeddings exist in \(V[G_\kappa][S][g]\).

	\[\begin{tikzcd}
	V[G_\kappa] \arrow[rr, "j"] \arrow[rdd, "i", swap] && M[G_\kappa][S][g][\Gtail]\\ \\ & N[G_\kappa][S'][g'][\Gtail'] \arrow[ruu, "k", swap]
	\end{tikzcd}\]	
	
	As the final act of forcing, let \(C\) be \(\C(S)^{V[G_\kappa][S]}\)-generic over \(V[G_\kappa][S][g]\).
	We claim that \(V[G_\kappa][S][g\times C]\) is our desired final extension. Recall that
	\cref{lemma:iterationResolvesToCohen}
	tells us that we can also write this extension as \(V[G_\kappa][H^0\times H^2]\) for some generic \(H^0\subseteq \Add(\kappa,\kappa^{++})^{V[G_\kappa]}\)
	and \(H^2\subseteq \Add(\kappa^{++},1)^{V[G_\kappa]}\). We will work from now on
	in this final model, using this alternative representation, and try to lift
	the embedding \(j\).
	
	By \cref{lemma:TermForcingCohen} we know that \(\Term(\P_\kappa,\Add(\kappa,\kappa^{++}))\) is forcing equivalent to
	\(\Add(\kappa,\kappa^{++})\) in \(V\). It follows from this by elementarity that the poset \(\Term(i(\P_\kappa), i(\Add(\kappa,\kappa^{++})))\)
	is equivalent to \(i(\Add(\kappa,\kappa^{++}))\) in \(N\). Now we return to an assumption we made at the
	start of the proof. 
	Since \(V\) has an \(i(\Add(\kappa,\kappa^{++}))\)-generic over \(N\), we can use
	this equivalence to also find a \(\Term(i(\P_\kappa), i(\Add(\kappa,\kappa^{++})))\)-generic over \(N\). Using \cref{fact:TermForcing}, we can combine
	this term forcing generic with the \(i(\P_\kappa)\)-generic \(G_\kappa*S'*g'*\Gtail'\)
	to extract an
	\(i(\Add(\kappa,\kappa^{++})^{V[G_\kappa]})\)-generic \(K'\)
	over \(N[G_\kappa][S'][g'][\Gtail']\) in \(V[G_\kappa][g']\). Since the forcing
	\(i(\Add(\kappa,\kappa^{++})^{V[G_\kappa]})\) is
	\(<i(\kappa)\)-distributive in \(N[G_\kappa][S'][g'][\Gtail']\), 
	\cref{fact:Transfer} again tells us that 
	the pointwise image \(k[K']\) generates a
	 \(j(\Add(\kappa,\kappa^{++})^{V[G_\kappa]})\)-generic filter \(\tilde{K}^0\) over \(M[G_\kappa][S][g][\Gtail]\).
	It is not necessarily the case that \(j[H^0]\subseteq \tilde{K}^0\), but we can surgically\footnote{See~\cite[Theorem 25.1]{Cummings:Handbook} or~\cite[Theorem 1, Second step]{Cummings1993:StrongUltrapowersLongCoreModels} for fairly detailed examples
	of this concrete use of the surgery method.} alter \(\tilde{K}^0\) 
	to obtain another \(j(\Add(\kappa,\kappa^{++})^{V[G_\kappa]})\)-generic \(K^0\) over
	\(M[G_\kappa][S][g][\Gtail]\) for which this will be the case, and we are able to lift \(j\) to
	\[j\colon V[G_\kappa][H^0]\to M[G_\kappa][S][g][\Gtail][K^0]\, .\]

	We can now forget about the maps \(i\) and \(k\) and focus solely on \(j\). To complete the lift,
	observe that \(\Add(\kappa^{++},1)^{V[G_\kappa]}\) remains \(\leq\kappa^+\)-distributive in
	\(V[G_\kappa][H^0]\) by Easton's lemma, 
	and so \cref{fact:Transfer} implies that the filter \(j[H^2]\) generates a generic \(K^2\) over
	\(M[G_\kappa][S][g][\Gtail][K^0]\), which gives us our final lift
	\[j\colon V[G_\kappa][H^0\times H^2]\to M[G_\kappa][S][g][\Gtail][K^0\times K^2]\, .\]
	Since \(j\) was originally a \((\kappa,\kappa^{++})\)-extender embedding, the same remains true for the lifted embedding, by \cref{fact:ExtenderLift}. 
	Since we clearly have
	\(2^\kappa=\kappa^{++}\) in the final model, \cref{lemma:extenderLiftsToUltrapower} tells us that the
	lift \(j\) is the ultrapower by a normal measure.
	
	\begin{lemma}
		\label{lemma:liftWitnessesCP}
		The embedding \(j\) witnesses \(\CP(\kappa,\kappa^+)\) in \(V[G_\kappa][H^0][H^2]\).
	\end{lemma}
	
	\begin{subproof}
		Let us write \(M^*=M[G_\kappa][S][g][\Gtail][K^0][K^2]\).
		We need to show that every subset \(x\) of \(\kappa^+\) in \(V[G_\kappa][S][C][g]\) appears in \(M^*\).
		To that end, we will first show that \(x\) is already in \(V[G_\kappa][S][g]\). This follows from \cref{lemma:ccDistributive}:
		the forcing to add \(C\) is \(<\kappa^{++}\)-distributive in \(V[G_\kappa][S]\), and \(\A(\kappa,\kappa^{++},\vec{Y})^{V[G_\kappa][S]}\)
		is \(\kappa^+\)-cc (and therefore trivially \(\kappa^{++}\)-cc) in \(V[G_\kappa][S][C]\),
		since it is equivalent to \(\Add(\kappa,\kappa^{++})\) in that model, as we explained in the proof of \cref{lemma:iterationResolvesToCohen}. Moreover, because of the distributivity of the forcing to add \(C\), the poset \(\A(\kappa,\kappa^{++},\vec{Y})^{V[G_\kappa][S]}\) is \(\kappa^{++}\)-cc in the model \(V[G_\kappa][S]\) as well.
		\cref{lemma:ccDistributive} then implies that the forcing to add \(C\)
		to \(V[G_\kappa][S][g]\) could not have added \(x\), and so \(x\) is already in that model.
		
		We next show that \(x\) has a name in \(M[G_\kappa]\). To start with,
		let \(\sigma\in V[G_\kappa][S]\) be a nice \(\A(\kappa,\kappa^{++},\vec{Y})\)-name for \(x\). Observe that
		\(\A(\kappa,\kappa^{++},\vec{Y})\) is actually a subset of \(H_{\kappa^{++}}^{V[G_\kappa]}\) (even though it is not an element of \(V[G_\kappa]\)), so the name \(\sigma\)
		is as well. Moreover, since \(\A(\kappa,\kappa^{++},\vec{Y})\) is
		\(\kappa^+\)-cc, \(\sigma\) has size \(\kappa^+\). But as a \(\kappa^+\)-sized subset of \(V[G_\kappa]\), the name \(\sigma\) could not
		have been added by the \(\leq\kappa^+\)-distributive forcing to add \(S\),
		and we conclude that \(\sigma\in H_{\kappa^{++}}^{V[G_\kappa]}\).
		Now, since \(\P_\kappa\) is \(\kappa\)-cc and \(\power(\power(\kappa))\in M\),
		we know that \(H_{\kappa^{++}}^{M[G_\kappa]}= H_{\kappa^{++}}^{V[G_\kappa]}\), so the name \(\sigma\) also appears in \(M[G_\kappa]\). 
		
		It follows that we can interpret the name \(\sigma\) by the generic filter
		\(g\) in \(M[G_\kappa][S][g]\) to find the set \(x\) in that model. Finally, we can conclude that \(M[G_\kappa][S][g]\) contains all the
		subsets of \(\kappa^+\) from \(V[G_\kappa][S][C][g]\), and so \(M^*\supseteq M[G_\kappa][S][g]\) does as well.
	\end{subproof}
		
	We have shown that \(\CP(\kappa,\kappa^+)\) holds in \(V[G_\kappa][H^0\times H^2]\). To finish the proof
	we also need to see that \(\kappa\) is the least measurable cardinal in that model. This
	follows easily from the way we designed the forcing \(\P_\kappa\). If \(\gamma<\kappa\)
	were measurable in \(V[G_\kappa][H^0\times H^2]\), it must definitely be inaccessible in \(V\).
	It follows that we did some nontrivial forcing at stage \(\gamma\) in the iteration 
	\(\P_\kappa\) and \cref{lemma:ApterShelahForcesNonmeasurable} implies that
	after the stage \(\gamma\) forcing \(\gamma\) is not measurable. The remaining forcing to get from
	that model to the model \(V[G_\kappa][H^0\times H^2]\) is at least \(\leq {2^\gamma}\)-strategically closed, which means
	that it could not have possibly added any measures on \(\gamma\). We can therefore conclude
	that \(\gamma\) remains nonmeasurable in \(V[G_\kappa][H^0\times H^2]\).
\end{proof}

The iteration we used is essentially the one described 
in~\cite[Section 4]{ApterShelah1997:MenasResultBestPossible}. It follows from the results outlined
there that, had we additionally assumed in \cref{thm:mainThmCP} that \(\kappa\) were 
\(\kappa^+\)-supercompact, this would remain true in the resulting extension.

\begin{corollary}
	If GCH holds and \(\kappa\) is \(\kappa^+\)-supercompact, 
	then there is a forcing extension in
	which \(\CP(\kappa,\kappa^+)\) holds, and \(\kappa\) is \(\kappa^+\)-supercompact and the least measurable. 
\end{corollary}

By starting with a stronger large cardinal hypothesis and modifying the forcing iteration
appropriately, we can push up the value of \(2^\kappa\) beyond just \(\kappa^{++}\) and capture
even more powersets.
In order to state the results as simply as possible, we make the following definition to add some convenient stages to the hierarchy of strong cardinals.

\begin{definition}
	If \(X\) is a set, a cardinal \(\kappa\) is called \emph{\(X\)-strong} if there
	is an elementary embedding \(j\colon V\to M\) with critical point \(\kappa\) and
	\(M\) a transitive inner model with \(X\in M\).
\end{definition}

\begin{theorem}
	\label{thm:CPwithLargeContinuum}
	Assume GCH holds and suppose that \(\kappa\) is \(H_\lambda\)-strong for some
	regular cardinal \(\lambda\geq\kappa^{++}\) which is \emph{not} the successor of a cardinal of cofinality less than \(\kappa\). Then there is a forcing extension in which
	\(\kappa\) is the least measurable cardinal, \(2^\kappa=\lambda\), and \(\CP(\kappa,<\lambda)\)
	holds (meaning that a single normal measure on \(\kappa\) captures every \(\power(\mu)\) for
	\(\mu<\lambda\)).
\end{theorem}

Note that this is a strict improvement over \cref{thm:mainThmCP} (we can pick \(\lambda=\kappa^{++}\) in the present theorem to recover the previous one). In particular, the value of \(2^\kappa=\lambda\) is optimal in the presence of \(\CP(\kappa,<\lambda)\), since this capturing property implies \(\CP(\kappa,\lambda')\) for each \(\lambda'<\lambda\), and those in turn imply that \(2^\kappa>\lambda'\), as we remarked earlier.

\begin{proof}
	The argument is much like the proof of \cref{thm:mainThmCP}, with a handful of changes: we will modify the forcing used slightly, and, more importantly, instead of preparing the model as
	in~\cite{Cummings1993:StrongUltrapowersLongCoreModels}, we use a preparation due to the second author. The different preparatory forcing also leads to the additional hypothesis on \(\lambda\) and, even assuming GCH in the ground model, will require some cardinal arithmetic calculations in order to be able to apply \cref{lemma:ApterShelahForcesNonmeasurable}. 
	
	As mentioned, we first use~\cite[Corollary 2.7]{Honzik2018:IndestructibilityForHypermeasurable} to pass to a forcing extension \(V^*\) of \(V\) in which the following hold:
		
	\begin{enumerate}
		\item \(2^\kappa=\kappa^+\) and \(2^{\kappa^+}=\lambda\).
		\item \(\kappa\) is \(H_\lambda\)-strong and this is witnessed by a \((\kappa,\lambda)\)-extender
		embedding \(j\colon V^*\to M\); moreover, \(M\) is closed under \(\kappa\)-sequences.
		\item\label{i:laver} There is a function \(\ell\colon \kappa\to\kappa\) such that \(\ell(\gamma)>\gamma^{++}\) is regular and not the successor of a cardinal of cofinality less than \(\gamma\) for all inaccessible \(\gamma<\kappa\), and \(j(\ell)(\kappa)=\lambda\).
		\item There is in \(V^*\) an \(M\)-generic filter for the poset \(j(\Add(\kappa,\lambda))\).
	\end{enumerate}

	The hypothesis on the values of \(\ell(\gamma)\) in (\ref{i:laver}) will allow us to conclude some cardinal arithmetic facts in \(V^*\), as we will describe momentarily. But first, let us briefly sketch the key parts of the preparatory forcing (details can be found in~\cite[Section 2.2]{Honzik2018:IndestructibilityForHypermeasurable}). First, we force, if necessary, with Woodin's fast function forcing to add a function \(\ell\) (see~\cite[Section 2.1]{Honzik2018:IndestructibilityForHypermeasurable}) which satisfies \(j(\ell)(\kappa)=\lambda\) for some embedding \(j\) witnessing the \(H_\lambda\)-strongness of \(\kappa\). This \(\ell\) may be assumed to have the properties described in the previous paragraph.
	
	The following forcing \(P^1\) is the Easton supported product of a collection of
	\(\leq\alpha\)-closed forcings \(Q^1_\alpha\), where \(\alpha\) runs through the set \(\mathcal{M}^+=\mathcal{M}\cup\{\kappa\}\), where \(\mathcal{M}\subseteq\kappa\) is the set of measurable closure points of the function \(\ell\). Each \(Q^1_\alpha\) is a lottery sum of forcing notions \(Q\) which are very close to being equal to \(\Add(\alpha^+,\ell(\alpha))\); more precisely, they are equal to \(i_W(\Add(\alpha,\ell(\alpha)))\), where \(i_W\) is the ultrapower embedding derived from some normal measure \(W\) on \(\alpha\) (see~\cite[Section 3.1]{Honzik2018:IndestructibilityForHypermeasurable} for more details on the connection between the \(Q\) and Cohen forcing at \(\alpha^+\)). In any case, the forcing notions constituting \(Q^1_\alpha\) live morally speaking on successors of cardinals in \(\mathcal{M}^+\), so the product-style definition of \(P^1\) is more natural (the Easton support \emph{iteration} is usually indicated when nontrivial forcing is done on a stationary set below \(\kappa\)). Additionally, and equally important, the product-style definition allows us to deal first with \(Q^1_\kappa\) and only later with the rest of \(P^1\), using the mutual genericity of the respective forcing notions. See~\cite[Lemma 2.3, Lemma 2.4]{Honzik2018:IndestructibilityForHypermeasurable} for more details.
	

	\begin{lemma}
		\label{lemma:cardinalArithmetic}
		In \(V^*\), if \(\gamma<\kappa\) is Mahlo and \(\theta<\ell(\gamma)\) then
		\(\theta^{<\gamma}<\ell(\gamma)\).
	\end{lemma}
	\begin{subproof}[Proof sketch]
		It will suffice to show that \(\theta^{\gamma'}<\ell(\gamma)\) for any Mahlo \(\gamma\), any \(\theta\) satisfying \(\gamma<\theta<\ell(\gamma)\), and any inaccessible \(\gamma'<\gamma\).
		Given such a \(\gamma'\), we can split up the product \(P^1=P^1_{<\gamma'}\times P^1_{\geq\gamma'}\) by grouping coordinates with indices less than \(\gamma'\) separately, and those with indices greater than or equal to \(\gamma'\) separately. The forcing \(P^1_{\geq\gamma'}\) is \(\leq\gamma'\)-closed, so it will not affect the value of \(\theta^{\gamma'}\). Let us focus on \(P^1_{<\gamma'}\).
		
		Since \(\gamma'\) is inaccessible in \(V^*\), we know that it cannot belong to any interval \((\alpha,\ell(\alpha)]\), since the forcing
		\(P^1\) forced \(2^{\alpha^+}=\ell(\alpha)\). In other words, we know
		that \(\ell(\beta)<\gamma'\) for all \(\beta<\gamma'\) in \(\mathcal{M}^+\). It follows that \(P^1_{<\gamma'}\) has size at most
		\(\gamma'\). From here, a simple calculation shows that there are at most
		\(\theta^{\gamma'}\) many nice \(P^1_{<\gamma'}\)-names for functions
		\(\gamma'\to\theta\). Since GCH holds in \(V\), we get \(\theta^{\gamma'}=\theta\) if \(\cf(\theta)>\gamma'\) and \(\theta^{\gamma'}=\theta^+\) if \(\cf(\theta)\leq\gamma'\). This means that, in \(V^*\), there are either \(\theta\) or \(\theta^+\) many functions \(\gamma'\to\theta\), depending on the cofinality of \(\theta\). Now recall that we wish to see that \(\theta^{\gamma'}<\ell(\gamma)\) in \(V^*\). We already know that \(\theta<\ell(\gamma)\), so the required inequality is immediate in the case that \(\cf(\theta)>\gamma'\). In the other case, when \(\cf(\theta)\leq\gamma'\), we recall that we assumed that \(\ell(\gamma)\) was not the successor of a cardinal of cofinality less than \(\gamma\). Since \(\cf(\theta)\leq\gamma'<\gamma\), it cannot be that \(\ell(\gamma)\) is equal to \(\theta^+\), so \(\theta^+<\ell(\gamma)\), as required.
%
%
	\end{subproof}

	A similar argument also shows that if \(\theta<\lambda\) then, in \(V^*\),
	we have \(\theta^{<\kappa}<\lambda\) for all \(\theta<\lambda\).
	
	The point of these calculations is to conclude that starting from \(V^*\),
	we can apply \cref{lemma:ApterShelahForcesNonmeasurable} to successively destroy the measurability of all cardinals below \(\kappa\), even without assuming GCH.

	Let us now move on from the preparatory forcing. In the interest of simpler notation, we will write just \(V\) instead of \(V^*\), and assume that all the properties enumerated above hold in \(V\).
	The initial iteration \(\P_\kappa\) will now be an Easton support iteration
	which forces at inaccessible cardinals \(\gamma<\kappa\) with the forcing \(\S_{\ell(\gamma)} * \A(\gamma,\ell(\gamma),\vec{X})\), with \(\vec{X}\) being
	an appropriate \(\clubsuit\)-sequence,
	\emph{provided that \(\gamma\) is inaccessible in \(V^{\P_\gamma}\)}.
	
	Since \(\P_\kappa\) factors at each inaccessible \(\gamma<\kappa\) into a two-step iteration of a small forcing (definitely of size less than \(\kappa\)) and a \(<\gamma\)-strategically closed forcing, we can readily see that \(\kappa\) remains a strong limit cardinal after forcing with \(\P_\kappa\).
	Moreover, since \(\P_\kappa\) is \(\kappa\)-cc, \(\kappa\) will remain inaccessible and there will be nontrivial forcing at stage \(\kappa\) of the iteration \(j(\P_\kappa)\), so we can write
	\[j(\P_\kappa)= \P_\kappa * \S_\lambda * \A(\kappa,\lambda,\vec{Y})* \Ptail\, .\]
	The full forcing that will give us the theorem is then
	\[\P= \P_\kappa * \S_\lambda * (\A(\kappa,\lambda,\vec{Y}) \times \C(\dot{S}))\, .\]
	The argument now proceeds very much like the proof of \cref{thm:mainThmCP},
	but with some simplifications due to the difference between the preparations
	from~\cite{Cummings1993:StrongUltrapowersLongCoreModels} and~\cite{Honzik2018:IndestructibilityForHypermeasurable}. We sketch the
	argument here, referring to the previously given proof and noting the main
	differences.
	
	Let \(G_\kappa*S*g\) be \(\P_\kappa*\S_\lambda*\A(\kappa,\lambda,\vec{Y})\)-generic
	over \(V\). We wish to lift the embedding \(j\colon V\to M\) to the extension
	\(V[G_\kappa]\) in the model \(V[G_\kappa][S][g]\). Given the factorization of \(j(\P_\kappa)\) above, we need to find a \(\Ptail\)-generic over \(M[G_\kappa][S][g]\). Previously we worked with the embeddings \(i\) and \(k\),
	but now we will be able to do without.\footnote{In fact, we could have employed the
	methods of~\cite{Honzik2018:IndestructibilityForHypermeasurable} even in the previous
	theorem, but we decided to give more details for the specific case \(2^\kappa=\kappa^{++}\). 
	}
	
	Consider any dense open subset \(D\)
	of \(\Ptail\) in \(M[G_\kappa][S][g]\). Since \(j\) was a \((\kappa,\lambda)\)-extender
	embedding, this \(D\) has the form \(D=j(f)(\alpha)^{G_\kappa*S*g}\) for some \(f\colon\kappa\to V_\kappa\) and some \(\alpha<\lambda\). 
	For each fixed \(f\) like this, the set
	\(\{j(f)(\alpha)^{G_\kappa*S*g}; \alpha<\lambda\}\) is an element of \(M[G_\kappa][S][g]\), since \(j(f)\in M[G_\kappa][S][g]\) and this model can evaluate this function and the resulting names and collect them together.
	Let \(\mathcal{D}_f\) be the subset of those elements of the form \(j(f)(\alpha)^{G_\kappa*S*g}\) that are dense open subsets of \(\Ptail\).
	Then \(\mathcal{D}_f\) is also an element of \(M[G_\kappa][S][g]\), since \(\Ptail\in M[G_\kappa][S][g]\).
%
	The set \(\mathcal{D}_f\) has size (at most) \(\lambda\) in \(M[G_\kappa][S][g]\). Since
	the first stage of forcing in \(\Ptail\) occurs beyond \(\lambda\), the forcing
	\(\Ptail\) is \(\leq\lambda\)-strategically closed in \(M[G_\kappa][S][g]\). This means that \(\bigcap\mathcal{D}_f\in M[G_\kappa][S][g]\) is a dense open subset of \(\Ptail\), and is also a subset of \(D\).
	
	Finally, observe that there are \(2^\kappa=\kappa^+\) many functions \(f\) (counted in \(V\)), and therefore only \(\kappa^+\) many dense sets \(\bigcap\mathcal{D}_f\).
	Since the forcing \(\P_\kappa*\S_\lambda*\A(\kappa,\lambda,\vec{Y})\) is composed
	of a \(\kappa^+\)-cc part, a \(<\lambda\)-distributive part, and another
	\(\kappa^+\)-cc part, applying \cref{fact:ClosurePreservation} twice allows
	us to conclude that \(M[G_\kappa][S][g]\) is closed under \(\kappa\)-sequences
	in \(V[G_\kappa][S][g]\). It follows that \(\Ptail\) remains \(\prec \kappa^+\)-strategically closed in \(V[G_\kappa][S][g]\), which will allow us to
	line up and meet all the dense sets \(\bigcap\mathcal{D}_f\) in turn, and so build
	a generic \(\Gtail\) for \(\Ptail\). This allows us to lift the embedding \(j\)
	to
	\[j\colon V[G_\kappa]\to M[G_\kappa][S][g][\Gtail]\]
	in \(V[G_\kappa][S][g]\).
	
	For the final step of the lift, we use \cref{lemma:iterationResolvesToCohen}
	to see \(\P\) as the iteration \(\P_\kappa*(\Add(\kappa,\lambda)\times\Add(\lambda,1))\).
	The lift through the forcing \(\Add(\kappa,\lambda)\) proceeds as in the proof of
	\cref{thm:mainThmCP}, except that we deal directly with the embedding \(j\)
	instead of passing through \(i\) and \(k\) as before. We apply 
	\cref{lemma:TermForcingCohen} to \(j(\Term(\P_\kappa,\Add(\kappa,\lambda)))\)
	in \(M\) and use our starting assumption that we have an \(M\)-generic for that poset
	in \(V\); a surgery argument like the one we alluded to before allows us to build
	a suitable \(j(\Add(\kappa,\lambda))\)-generic \(K^0\) over \(M[G_\kappa][S][g][\Gtail]\) and lift \(j\) to
	\[j\colon V[G_\kappa][H^0]\to M[G_\kappa][S][g][\Gtail][K^0]\,.\]
	The lift through the final forcing \(\Add(\lambda,1)^{V[G_\kappa]}\) is handled
	exactly as in the proof of \cref{thm:mainThmCP}. 
	
	It remains for us to see that this lifted
	embedding witnesses \(\CP(\kappa,<\lambda)\) in the final model and that \(\kappa\)
	is the least measurable cardinal there.
	
	We need to show that every bounded subset \(x\) of \(\lambda\) in \(V[G_\kappa][S][g][C]\) appears in the target model of the lifted embedding \(j\). This works almost exactly as in \cref{lemma:liftWitnessesCP}. We first use \cref{lemma:ccDistributive} to show that \(x\) is already in \(V[G_\kappa][S][g]\). Then we argue that \(x\) has a name in \(M[G_\kappa]\). This is because we can find a nice
	\(\A(\kappa,\lambda,\vec{Y})\)-name \(\sigma\) for \(x\) in \(V[G_\kappa][S]\), of size less than \(\lambda\), that is a subset of \(H_{\lambda}^{V[G_\kappa]}\). Since the forcing to add \(S\) is \(<\lambda\)-distributive, this name could not have been added by it, so
	\(\sigma\in H_\lambda^{V[G_\kappa]}\). But since \(\P_\kappa\) is \(\kappa\)-cc and \(H_\lambda^V=H_\lambda^M\) (as \(j\) witnesses the \(H_\lambda\)-strongness of \(\kappa\)), it follows that \(H_\lambda^{V[G_\kappa]}=H_\lambda^{M[G_\kappa]}\). As in the previous proof, this means that we can interpret \(\sigma\) in \(M[G_\kappa][S][g]\)
	to find \(x\) in that model, as well as in the target model of \(j\). So the lifted embedding \(j\) really does witness \(\CP(\kappa,<\lambda)\).
	
	To see that \(\kappa\) is the least measurable cardinal in the final model, we simply inspect our construction of the forcing \(\P_\kappa\). If \(\gamma<\kappa\) were measurable in the final model, it must necessarily be Mahlo in the intermediate extension \(V[G_\gamma]\), and so some nontrivial forcing must have occurred in the next step. It follows from \cref{lemma:cardinalArithmetic} that, over this model, the next step of forcing with \(\A(\gamma,\ell(\gamma),\vec{X})\) destroyed the measurability of \(\gamma\), and the remainder of the forcing possesses too much closure to ever recover this measurability.
\end{proof}

Conversely, we can extend Cummings' argument to show that the large cardinal hypothesis we used
above is optimal.

\begin{theorem}
	\label{thm:consistencyStrengthOfCP}
	Suppose that \(\CP(\kappa,<\lambda)\) holds for some regular cardinal \(\lambda\geq \kappa^{++}\).
	Then \(\kappa\) is \(H_\lambda\)-strong in an inner model. Moreover, this inner model
	satisfies GCH, and so \(\kappa\) is \((\kappa+\alpha)\)-strong there, where
	\(\lambda=\kappa^{+\alpha}\).
\end{theorem}

\begin{proof}
	This is essentially standard. Suppose that \(j\colon V\to M\) is an ultrapower embedding by a normal measure witnessing \(\CP(\kappa,<\lambda)\); it follows that
	\(H_\lambda\in M\).\footnote{Recall that \(\CP(\kappa,<\lambda)\) implies \(2^{\kappa}\geq \lambda\), so \(H_\lambda\) being in \(M\) is weaker than \(V_{\kappa+\alpha}\) being in \(M\), where \(\lambda=\kappa^{+\alpha}\). In particular, \(\kappa\) might not be \((\kappa+\alpha)\)-strong in \(V\).} 
	We assume that there is no
	inner model with a strong cardinal and let \(K\) be the core model with the (nonoverlapping) extender sequence \(\vec{E}\). It follows
	that \(j\rest K\) is the result of a normal iteration of \(\vec{E}\) and, since the critical point
	of \(j\) is \(\kappa\), the first extender applied in this iteration must have index
	\((\kappa,\eta)\) for some \(\eta\). Since \(\vec{E}\) is coherent, the sequence \(j(\vec{E})\)
	has no extenders with indices \((\kappa,\beta)\) for \(\beta\geq\eta\). But since
	\(M\) captured all of \(H_\lambda\), we must have \(K\rest\lambda=K^M\rest\lambda\), and
	so \(\vec{E}\) and \(j(\vec{E})\) must agree up to \(\lambda\). It follows that \(\eta\geq\lambda\)
	and so \(o(\kappa)\geq\lambda+1\) (and \(\kappa\) is \(H_\lambda\)-strong) in \(K\).
	
	Since \(K\) satisfies GCH, \(V_{\kappa+\alpha}^K\) is a transitive set of size \(\kappa^{+\alpha}=\lambda\) there. It follows that the transitive closure
	of each element of \(V_{\kappa+\alpha}^K\) has size strictly less than \(\lambda\),
	so these elements appear in the codomain of the embedding witnessing the \(H_\lambda\)-strongness of \(\kappa\) in \(K\).
\end{proof}

The preparation from~\cite{Honzik2018:IndestructibilityForHypermeasurable} works even for singular \(\lambda\)
of cofinality strictly above \(\kappa\) (if the cofinality of \(\lambda\) is equal to \(\kappa^{+}\),
we get \(2^{\kappa^+}=\lambda^+\) in (1) above). It is unclear, however, whether \cref{thm:CPwithLargeContinuum}
can allow for this weaker hypothesis (in particular, \cref{lemma:iterationResolvesToCohen} 
seems to rely crucially on the second parameter in the forcing \(\A\) being regular).

\begin{question}
	Can \cref{thm:CPwithLargeContinuum} be improved to allow for arbitrary \(\lambda\)
	of cofinality strictly above \(\kappa\)?
\end{question}

Another question raised by \cref{thm:CPwithLargeContinuum} is whether \(\CP(\kappa,\lambda)\)
can fail for the first time at some \(\kappa^+<\lambda<2^\kappa\). The following theorem shows that
the answer is yes.
%

\begin{theorem}
	\label{thm:intermediateFailureOfCP}
	Suppose that there is no inner model with a strong cardinal and let \(V=K\) be the core model. 
	Suppose that \(\kappa\) is \(H_{\kappa^{+3}}\)-strong. 
	Then there is a forcing extension in which \(\kappa\) is the least measurable cardinal, \(2^\kappa=\kappa^{+3}\), and \(\CP(\kappa,\kappa^+)\) holds while \(\LCP(\kappa,\kappa^{++})\) fails.
\end{theorem}

\begin{proof}
	We will use the same forcing as in the proof of \cref{thm:CPwithLargeContinuum}, letting \(\lambda=\kappa^{+3}\) (note that the core model satisfies GCH, so the hypotheses of that theorem are satisfied). 
	That is, we shall force with \[\R=P^1*\P_\kappa*\S_{\kappa^{+3}}*(\A(\kappa,\kappa^{+3},\vec{Y})\times \C(\dot{S}))\,,\]
	using the notation from the proof of \cref{thm:CPwithLargeContinuum}.
	We already know that after forcing with this \(\R\) we obtain \(2^\kappa=\kappa^{+3}\) and \(\CP(\kappa,\kappa^{++})\), while \(\kappa\) becomes the least measurable cardinal.
	To obtain the desired extension, we shall force with \(\R\times \Add(\kappa^{++},1)^V\).
	
	\begin{lemma}
		\label{lemma:Easton}
		The forcing \(\Add(\kappa^{++},1)^V\) remains \(\leq\kappa^+\)-distributive in \(V^\R\).
	\end{lemma}
	\begin{subproof}
		This is essentially a version of Easton's lemma.
		Let us write \(\R'=P^1*\P_\kappa\).
		We can rewrite \(\R\) as
		\[\R'*(\Add(\kappa,\kappa^{+3})\times \Add(\kappa^{+3},1))\,,\]
		using \cref{lemma:iterationResolvesToCohen}.
		It follows from~\cite[Lemma 2.3]{Honzik2018:IndestructibilityForHypermeasurable} that \(\R'*\Add(\kappa,\kappa^{+3})\) is \(\kappa^{++}\)-cc, and, of course, \(\Add(\kappa^{+3},1)\) is \(\leq\kappa^{++}\)-closed in \(V^{\R'}\).
		
		Now let \(\dot{f}\) be an \(\R\times \Add(\kappa^{++},1)^V\)-name for a \(\kappa^+\)-sequence of ordinals. 
		For an ordinal \(\alpha<\kappa^+\), a condition \(s\in\Add(\kappa^{++},1)^V\), and an \(\R'\)-name \(\dot{q}\) for a condition in \(\Add(\kappa^{+3},1)^{V^{\R'}}\), say that the pair \((\dot{q},s)\) is \emph{\(\alpha\)-good} if there is a maximal antichain of conditions \((r,\dot{p})\in \R'*\Add(\kappa,\kappa^{+3})\) such that \((r,\dot{p},\dot{q},s)\) decides the value of \(\dot{f}(\alpha)\). 
		We will see that any condition in \(\R\times \Add(\kappa^{++},1)^V\) can be strengthened to one whose latter two coordinates are \(\alpha\)-good for all \(\alpha<\kappa^+\).
		
		Pick \(\dot{q}\) and \(s\) as above and consider the set of all \((r,\dot{p})\) such that \((r,\dot{p},\dot{q},s)\) decides \(\dot{f}(0)\). 
		This is an open set of conditions, so we may pick a maximal antichain \(W^0_0\) from this set, and it will remain an antichain as a subset of \(\R'*\Add(\kappa,\kappa^{+3})\). 
		If \(W^0_0\) is already a maximal antichain in \(\R'*\Add(\kappa,\kappa^{+3})\), then \((\dot{q},s)\) is 0-good.
		Otherwise we can find some \((r_0,\dot{p}_0)\) which is incompatible with every condition in \(W^0_0\). 
		We can also find some \((r_1,\dot{p}_1,\dot{q}',s_1)\leq (r_0,\dot{p}_0,\dot{q},s)\) which decides \(\dot{f}(0)\).
		Using a mixing argument, let \(\dot{q}_1\) be an \(\R'\)-name such that \(\R'\forces \dot{q}_1\leq \dot{q}\) and \(r_1\forces \dot{q}_1=\dot{q}'\).
		Now consider the set of all \((r,\dot{p})\) such that \((r,\dot{p},\dot{q}_1,s_1)\) decides \(\dot{f}(0)\). 
		This set includes \(W_0^0\) as well as \((r_1,\dot{p}_1)\), so we may again pick a maximal antichain \(W_1^0\supset W_0^0\) from it. 
		If \(W_1^0\) turns out to be maximal in \(\R'*\Add(\kappa,\kappa^{+3})\), then \((\dot{q}_1,s_1)\) is 0-good, and otherwise we can keep going.
		
		We continue recursively, constructing larger and larger antichains \(W_\eta^0\). At limit stages we take unions of the previously constructed antichains and use the \(\leq\kappa^+\)-closure of \(\Add(\kappa^{++},1)^V\)
		and the (forced) \(\leq\kappa^{++}\)-closure of \(\Add(\kappa^{+3},1)^{V^{\R'}}\) to find lower bounds for the sequences of conditions \(s_\eta\) and \(\dot{q}_\eta\). 
		The closure suffices to continue this construction for all \(\eta<\kappa^{++}\) (although notice that the degree of closure in the \(s\) component is too low to find a putative lower bound \(s_{\kappa^{++}}\)). 
		However, the construction must in fact stop at some stage before \(\kappa^{++}\), otherwise the union \(W^0=\bigcup_{\eta<\kappa^{++}}W^0_\eta\) would be an antichain in \(\R'*\Add(\kappa,\kappa^{+3})\) of size \(\kappa^{++}\), contradicting the chain condition of that poset.
		Once the construction stabilizes, we've reached a 0-good pair \((\dot{q}_\eta,s_\eta)\), as witnessed by the maximal antichain \(W^0_\eta\). 
		Notice that \(s_\eta\leq s\) and \(\forces_{\R'}\dot{q}_\eta\leq \dot{q}\).
		
		Repeating the same argument for all \(\alpha<\kappa^+\), we see that, starting with any \(\dot{q}\) and \(s\), we can find \(\forces_{\R'}\dot{q}'\leq \dot{q}\) and \(s'\leq s\) such that \((\dot{q}',s')\) is \(\alpha\)-good. Since any pair stronger than an \(\alpha\)-good pair is itself \(\alpha\)-good, we can use closure in both coordinates one last time to find, below any \((\dot{q},s)\), a pair which is \(\alpha\)-good for all \(\alpha<\kappa^+\).
		
		Finally, let \(G\times H\subseteq \R\times \Add(\kappa^{++},1)^V\) be generic.
		By the density property just described, we can find a condition \((r,\dot{p},\dot{q},s)\in G\times H\) whose latter two coordinates are \(\alpha\)-good for all \(\alpha<\kappa^+\). 
		But given such a condition, we can find \(f(\alpha)\) by consulting where the generic \(G\) (or even \(G\rest \R'*\Add(\kappa,\kappa^{+3})\)) meets the maximal antichain witnessing the \(\alpha\)-goodness of \((\dot{q},s)\).
		Therefore we can find \(f\in V[G]\), which is what we needed to show.
	\end{subproof}

	The poset \(\Add(\kappa^{++},1)^V\) has size \(\kappa^{++}\) in \(V\) (since we have GCH), so it remains \(\kappa^{+3}\)-cc in the extension by \(\R\). Since \cref{lemma:Easton} shows that this poset is also \(\leq\kappa^+\)-distributive in the extension by \(\R\), it follows that cardinals are preserved
	to the final extension by \(\R\times \Add(\kappa^{++},1)^V\) and that \(\kappa\) remains the least measurable cardinal.
	Moreover, the ultrapower embedding witnessing \(\CP(\kappa,\kappa^{++})\) lifts to the extension by \cref{fact:Transfer} and, since the extension by \(\Add(\kappa^{++},1)^V\) does not add any subsets of \(\kappa^+\), the lifted embedding still witnesses \(\CP(\kappa,\kappa^+)\).
	
	\begin{lemma}
		\(\LCP(\kappa,\kappa^{++})\) fails in the final extension.
	\end{lemma}
	\begin{subproof}
		Recall that we are forcing over \(V=K\) using the poset
		\[\left(P^1*\P_\kappa*\S_{\kappa^{+3}}*(\A(\kappa,\kappa^{+3},\vec{Y})\times \C(\dot{S}))\right)\times \Add(\kappa^{++},1)^V\,,\]
		which we can rewrite in the form
		\[\left(P^1*\P_\kappa*(\Add(\kappa,\kappa^{+3})\times \Add(\kappa^{+3},1))\right)\times \Add(\kappa^{++},1)^V\,.\]
		Let \((G^1*G_\kappa*(H^0\times H^3))\times H^2\) be generic over \(V\) for this poset.
		Now suppose that \(\LCP(\kappa,\kappa^{++})\) holds in the final extension \(V[(G^1*G_\kappa*(H^0\times H^3))\times H^2]\).
		This means that there is a normal ultrapower
		\[j^*\colon V[G^1][G_\kappa][H^0\times H^3\times H^2] \to M^*[g^1][G^*][H^*]\]
		on \(\kappa\) which captures both \(H^2\subseteq \kappa^{++}\) and \(\power(\kappa^+)^V\) (which can be coded as a subset of \(\kappa^{++}\), since GCH holds in \(V\)). 
		Above, we intended \(g^1\) to be \(j^*(P^1)\)-generic over \(M^*\), \(G^*\) to be \(j^*(\P_\kappa)\)-generic over \(M^*[g_1]\), and \(H^*\) to be generic over \(M^*[g^1][G^*]\) for the remaining Cohen forcing. 
		Since no cardinals were collapsed between \(V\) and the final extension and we insisted that \(\power(\kappa^+)^V\in M^*[g^1][G^*][H^*]\), we can conclude that this model computes \(\kappa^{++}\) correctly.
		
		Since \(j^*\) is an elementary embedding of a generic extension of the core model, its restriction to \(j^*\colon V\to M^*\) is an iteration of the core model itself, and therefore \(M^*=V\cap M^*[g^1][G^*][H^*]\) is an inner model of \(V\).\footnote{This is the only place we use the fact that we started this construction in the core model.} 
		This implies that \(\power(\kappa^+)^V\in M^*\) (since we explicitly put this powerset into \(M^*[g^1][G^*][H^*]\)) and also that GCH holds in \(M^*\).
		
		Consider the forcing \(j^*(P^1)=P^1\rest\kappa \times P^{1*}_\kappa \times P^{1*}_{\text{tail}}\), where \(P^{1*}_\kappa\) is the factor of \(j^*(P^1)\) indexed at \(\kappa\). 
		Note that \(P^1\rest\kappa\) really is an initial segment of \(j^*(P^1)\), since we necessarily have \(p=j^*(p)\in j^*(P^1)\) for \(p\in P^1\rest\kappa\). 
		The forcing \(P^1\rest\kappa\) has size \(\kappa\), and it follows from the precise description of \(Q^1_\alpha\) in~\cite[Section 2.2]{Honzik2018:IndestructibilityForHypermeasurable} and \cref{lemma:GitikMerimovich} that \(P^{1*}_\kappa\) is equivalent, over \(M^*\), to \(\Add(\kappa^+,\kappa^{+3})^{M^*}\). Since \(M^*\) satisfies GCH, this poset is \(\kappa^{++}\)-Knaster in \(M^*\). As the product of two \(\kappa^{++}\)-Knaster posets, \(P^1\rest\kappa \times P^{1*}_\kappa\) is itself \(\kappa^{++}\)-Knaster in \(M^*\), which implies that its square is \(\kappa^{++}\)-cc in \(M^*\).
		
		Unger~\cite[Lemma 2.4]{Unger2013:AronszajnTreesAndSuccessorsOfSingulars} showed that any poset whose square was \(\lambda\)-cc for some regular \(\lambda\) has the \(\lambda\)-approximation property, which states that any set of ordinals in the extension, all of whose subsets of size less than \(\lambda\) are in the ground model, must itself be in the ground model. As a special case, such forcings cannot add fresh subsets of \(\lambda\) (recall that a set of ordinals is \emph{fresh} over a model if it is not in that model but all of its initial segments are).
		Applying this to our situation, we can see that forcing over \(M^*\) by \(P^1\rest \kappa\times P^{1*}_\kappa\) does not add any new fresh subsets of \(\kappa^{++}\). Of course, \(H^2\) is a fresh subset of \(\kappa^{++}\) over \(V\), and since \(V\) and \(M^*\) have the same bounded subsets of \(\kappa^{++}\), it is also fresh over \(M^*\). Therefore \(H^2\) is not added to \(M^*\) by \(P^1\rest\kappa \times P^{1*}_\kappa\). Moreover, the tail forcing \(P^{1*}_{\text{tail}}\) is \(\leq \kappa^{++}\)-closed over
		\(M^*\) by~\cite[Lemma 2.3]{Honzik2018:IndestructibilityForHypermeasurable}, and therefore does not add any subsets of \(\kappa^{++}\) to the extension of \(M^*\) by \(P^1\rest\kappa \times P^{1*}_\kappa\) by Easton's lemma. Hence, \(H^2\) does not appear in \(M^*[g^1]\).
		
		Now let us write \(G^*=G_\kappa*(S^**g^**\Gtail^*)\), where \(S^*\) is generic for \(\S_{\kappa^{+3}}\), and \(g^*\) is generic for \(\A(\kappa,\kappa^{+3},\vec{Y})\). As was the case above, \(G_\kappa\) really is an initial segment of \(G^*\), for the same reason. The forcing \(G_\kappa\) has size \(\kappa\) and therefore cannot add \(H^2\) to \(M^*[g^1]\) by another application of Unger's result. 
		On the other hand, the forcing to add \(S^*\) is \(\leq \kappa^{++}\)-distributive and therefore also cannot add \(H^2\).
		
		Note that \(M^*[g^1]\) satisfied \(2^\kappa=\kappa^{+}\) and \(2^{\kappa^+}=\kappa^{+3}\) by~\cite[Corollary 2.7]{Honzik2018:IndestructibilityForHypermeasurable}. This remains true after adding \(G_\kappa*S^*\) as well, and \cref{cor:ApterShelahKnaster} tells us that the forcing \(\A(\kappa,\kappa^{+3},\vec{Y})\) is \(\kappa^{++}\)-Knaster, which implies that its square is \(\kappa^{++}\)-cc. Applying Unger's result yet again, we see that \(H^2\) does not appear in \(M^*[g^1][S^**g^*]\). On the other hand, the tail forcing adding \(\Gtail^*\) is \(\leq\kappa^{++}\)-closed, and also does not add \(H^2\).
		
		Finally, observe that the forcing to add \(H^*\) is at least \(<j^*(\kappa)\)-distributive, so it definitely cannot add \(H^2\) to \(M^*[g^1][G^*]\). But this contradicts our original assumption that \(H^2\in M^*[g^1][G^*][H^*]\).
%
	\end{subproof}
	To summarize, we've obtained a forcing extension of \(V\) in which \(\kappa\) is the least measurable cardinal, \(\CP(\kappa,\kappa^+)\) holds, but \(\LCP(\kappa,\kappa^{++})\) fails, as required.
\end{proof}

One would expect that it should be possible to force \(2^\kappa=\kappa^{+3}\)
and \(\CP(\kappa,\kappa^+)\) starting from a large cardinal hypothesis weaker than an \(H_{\kappa^{+3}}\)-strong cardinal \(\kappa\); an \(H_{\kappa^{+2}}\)-strong and \(\kappa^{+3}\)-tall cardinal \(\kappa\) likely suffices 
(recall that \(\kappa\) is \emph{\(\lambda\)-tall} if there is an elementary embedding \(j\colon V\to M\) with critical point \(\kappa\) such that \(M\) is closed under \(\kappa\)-sequences and \(j(\kappa)>\lambda\); see~\cite{Hamkins2009:TallCardinals}).

\begin{question}
	What is the consistency strength of \(\CP(\kappa,\kappa^+)\) and \(2^\kappa=\kappa^{+3}\) holding at the least measurable cardinal \(\kappa\) but
	\(\LCP(\kappa,\kappa^{++})\) failing?
\end{question}

It is also unclear whether the anti-large cardinal hypothesis and use of the core model are crucial for the above result. The only use of that hypothesis comes when we wish to understand the nature of generic embeddings of the ground model. It is plausible that one could use Hamkins' results on elementary embeddings in generic extensions with the approximation and cover properties (see~\cite{Hamkins2003:ApproximationAndCoveringNoNewLargeCardinals}) to prove a more general result, but those theorems do not interact well with the product nature of our preparation \(P^1\). One potential approach would be to mimic the proof of \cref{thm:mainThmCP} but to weave more complicated forcing into the preparation.

\begin{question}
	Can one obtain a model as in \cref{thm:intermediateFailureOfCP} without starting from the core model?
\end{question}

At the end of the paper, let us give another example of the power of \cref{lemma:extenderLiftsToUltrapower}
in showing that \(\CP(\kappa,\kappa^+)\) holds in known forcing extensions. As we have seen,
\(\CP(\kappa,\kappa^+)\) does not have any implications for the outright size of \(\kappa\), since
it may consistently hold at the least measurable cardinal \(\kappa\). But one might try to measure its effects
slightly differently. While the capturing property says that there is a normal measure on \(\kappa\) which
is quite ``fat'', in the sense that it captures all subsets of \(\kappa^+\), perhaps \(\kappa\) must inevitably
also carry some, or many, ``thin'' measures which do not capture much at all. In other words,
perhaps \(\CP(\kappa,\kappa^+)\) has some implications about the number of normal measures on \(\kappa\).
A combination of \cref{lemma:extenderLiftsToUltrapower} and a theorem of Friedman and Magidor will show
us that this is not the case.

\begin{theorem}
	\label{thm:FriedmanMagidor}
	If \(V\) is the minimal extender model with a \((\kappa+2)\)-strong cardinal \(\kappa\) and \(\lambda\leq\kappa^{++}\) 
	is a cardinal, then there is a forcing extension in which \(\kappa\) carries exactly \(\lambda\) many normal measures
	and each of them witnesses \(\CP(\kappa,\kappa^+)\). In particular, it is consistent that \(\kappa\) has a unique
	normal measure and \(\CP(\kappa,\kappa^+)\) holds.
\end{theorem}

\begin{proof}
	The hard part of the proof was done by Friedman and Magidor~\cite[Theorem 19]{FriedmanMagidor2009:NumberNormalMeasures}, 
	who showed
	that, starting from a model \(V\) as in the hypothesis of this theorem, there is a forcing extension \(V[G]\) satisfying \(2^\kappa=\kappa^{++}\) 
	in which \(\kappa\) carries exactly \(\lambda\) many normal measures. 
	They also show that each of these normal measures is derived from a lift of the ground model extender embedding \(j\colon V\to M\)
	witnessing the \((\kappa+2)\)-strongness of \(\kappa\). However, \cref{lemma:extenderLiftsToUltrapower} implies
	that these lifts are themselves already ultrapowers by a normal measure on \(\kappa\). Finally, an analysis of their
	proof shows that the forcing used to obtain the model \(V[G]\) can be written as \(\P*\dot{\Q}\) where \(\P\subseteq H_{\kappa^{++}}\) is a \(\kappa^{++}\)-cc poset
	which is regularly embedded in \(j(\P)\), and \(\dot{\Q}\) is forced to be \(\leq\kappa^+\)-distributive. 
	It follows that every subset of \(\kappa^+\) in \(V[G]\) has a nice name in \(H_{\kappa^{++}}^V\in M\) and therefore
	appears in \(M[j(G)]\).
\end{proof}

It is unclear whether one can obtain similar results at the least measurable cardinal \(\kappa\).
It seems likely that, to do so, it would be necessary to adapt the forcing \(\A\)
to incorporate the Sacks forcing machinery that Friedman and Magidor used in their arguments.

\begin{question}
	Is it consistent that the least measurable cardinal \(\kappa\) carries a unique normal
	measure and \(\CP(\kappa,\kappa^+)\) holds?\footnote{In a personal communication, James Cummings has answered this question in the affirmative. He showed that in the Friedman--Magidor model described in \cref{thm:FriedmanMagidor} \(\kappa\) is the \emph{only} measurable cardinal.}
\end{question}

\bibliographystyle{amsplain}
\bibliography{bibbase}
\end{document}